\RequirePackage{fix-cm}
\documentclass[smallextended]{svjour3}       % onecolumn (second format)
\smartqed  % flush right qed marks, e.g. at end of proof
\usepackage{cite}
\usepackage{graphicx}
\usepackage{graphicx,amssymb,amsmath}
\usepackage{booktabs}
\usepackage{array}
\usepackage{color}
\usepackage{enumerate}
\usepackage{algorithm}
\usepackage{algorithmic}
\usepackage{float}
\usepackage{multirow}
\usepackage{rotating}
\usepackage{epstopdf}
\usepackage[misc]{ifsym}
\usepackage{bm}
\usepackage[colorlinks,linkcolor=blue,citecolor=black]{hyperref}
\begin{document}

\title{PAL-Hom method for QP and an application to LP%\thanks{Grants or other notes
}
\subtitle{}


\author{Guoqiang Wang \and Bo Yu}


\institute{Guoqiang Wang \at
School of Mathematical Sciences, Dalian University of Technology, Dalian, Liaoning 116024, P.
R. China             \\
              \email{wangguojim@mail.dlut.edu.cn}           %  \\
           \and
           Bo Yu({\Letter}) \at
School of Mathematical Sciences, Dalian University of Technology, Dalian, Liaoning 116024, P.
R. China           \\
           \email{yubo@dlut.edu.cn}
}

\date{Received: date / Accepted: date}

\maketitle

\begin{abstract}
In this paper, a proximal augmented Lagrangian homotopy (PAL-Hom) method for solving convex quadratic programming problems is  proposed. This method takes the proximal augmented Lagrangian method as the outer iteration. To solve the proximal augmented Lagrangian subproblems, a homotopy method is presented as the inner iteration. The homotopy method tracks the piecewise-linear solution path of a parametric quadratic programming problem whose start problem takes an approximate solution as its solution and the target problem is the subproblem to be solved. To improve the performance of the homotopy method, the accelerated proximal gradient method is used to obtain a fairly good approximate solution that implies a good prediction of the optimal active set. Moreover, a sorting technique for the Cholesky factor update  as well as an $\varepsilon$-relaxation technique for checking primal-dual feasibility and correcting the active sets are presented to improve the efficiency and robustness of the homotopy method. Simultaneously, a proximal-point-based AL-Hom method  which is shown to converge in finite number of steps, is applied to linear programming. Numerical experiments on randomly generated  problems and the problems from the CUTEr and  Netlib test collections, support vector machines (SVMs) and contact problems of elasticity demonstrate that PAL-Hom is faster than the active-set methods and the parametric active set methods and is competitive to the interior-point methods and the specialized algorithms designed for specific models (e.g., sequential minimal optimization (SMO) method for SVMs).

\end{abstract}

\section{Introduction}

In this paper, we consider the convex quadratic programming (QP) problem
\begin{equation}\label{scqp}
\begin{array}{l}
\min ~~ \frac{1}{2}x^{T}Qx+r^{T}x \\	
\rm{s.t.}~~~\emph{Ax}=\emph{b},\\
~~~~~~~\emph{x}\geq 0,
\end{array}
\end{equation}
where $Q$ is an $n\times n$ symmetric semipositive definite matrix, $A$ is an $m\times n$ matrix, $r$ is an $n$-dimensional vector and $b$ is an $m$-dimensional vector.

As a type of classical optimization problem, QP problems arise in many areas, e.g., finance \cite{cornuejols2006optimization}  and optimal control \cite{dostal2000solution, dostal2000duality}. In particular, the QP problem is a key issue in support vector machines (SVMs) \cite{CC01a}. Due to their broad applicability, QP problems have attracted an enormous amount of research  on developing  efficient algorithms. Among the typical methods for QP problems, interior-point methods (IPMs) \cite{karmarkar1984new, wright1997primal, mehrotra1992implementation, zhang1998solving} and active-set (AS) methods  \cite{fletcher1971general, fletcher2000stable, gould1991algorithm, forsgren2016primal, gill2015methods} are two important types of methods and have been implemented in many software packages, e.g., CPLEX, Gurobi, MATLAB, IPOPT\cite{wachter2006implementation}, QPOPT \cite{gill1995user}, and SNOPT\cite{gill2008user}.  AS methods are efficient for solving small- to medium-scale QPs, and generally, they can obtain high-precision solutions. Compared with AS methods, IPMs are shown to be more efficient for large-scale QP problems.

In addition to these two classical methods, the augmented Lagrangian method (ALM) is another well-known method and was proposed independently by Hestenes \cite{hestenes1969multiplier} and Powell \cite{powell1969method} for nonlinear programming with general constraints and simple bounds. Simultaneously, Powell provided a global convergence analysis that requires the exact minimization of the subproblems. However, Rockafellar \cite{rockafellar1976augmented}, Bertsekas \cite{bertsekas1999nonlinear}  and Conn et al. \cite{conn1991globally} showed that  the exact minimization is not necessary for convergence. Moreover, a superlinear convergence analysis of the ALM that exactly solves the subproblems and has a new Lagrangian multipliers updating formula was presented by Yuan \cite{ yuan2014analysis}. Besides the theory results, Conn presented an effective numerical implementation of the ALM that exactly asymptotically solved the subproblems in  LANCELOT software \cite{conn2013lancelot}.

Dostal et al. used the ALM \cite{dostal2003augmented} to solve QP problems  $(\ref{scqp})$, which follows Conn et al.  \cite{conn1991globally}  who used the ALM to solve general nonlinear constrained optimization problems.
The $k$-th iteration of the ALM for $(\ref{scqp})$ begins with a given $\lambda^{k}$ and obtains $(x^{k+1},~\lambda^{k+1})$ via
\begin{eqnarray}\label{alm subproblem}
x^{k+1} & =&\arg\min~\{\mathcal{L}_{\beta}(x,\lambda^{k})~|~x\geq0\},\\ \nonumber
\lambda^{k+1} & =&\lambda^{k}-\beta(Ax^{k+1}-b).
\end{eqnarray}
where
\begin{eqnarray}\nonumber
\mathcal{L}_{\beta}(x,\lambda)=\frac{1}{2}x^{T}Qx+r^{T}x-\lambda^{T}(Ax-b)+\frac{\beta}{2}\Arrowvert Ax-b \Arrowvert^{2}
\end{eqnarray}
is the augmented Lagrangian function of $(\ref{scqp})$ (omitting the $x\geq 0$ bounds). The performance of the ALM depends on the solving of the subproblems (\ref{alm subproblem}). Therefore, it is desired to design efficient algorithms for the augmented Lagrangian subproblems.

The parametric active-set (PAS) method is a type of AS method proposed by Ritter \cite{ritter1967method, ritter1981parametric} and Best \cite{best1982algorithm, best1996algorithm} for the parametric quadratic programming (PQP) problem
\begin{equation}\label{pqp0}
\begin{array}{l}
\min ~\{(r+tq)^{T}x+\frac{1}{2}x^{T}Q{x} |~g+tp\leq Bx\}
\end{array}
\end{equation}
where $p$ and $q$ are $n$-vectors, $g$ is an $m$-vector, and $B$ is an $m\times n$ matrix. Ferreau et al. \cite{ferreau2008online} applied the PAS method to the model predictive control problem by solving a sequence of PQP  problems with the PAS method. The PQP problems are constructed such that the starting solution of the current PQP problem occurs as the target solution of the previous PQP problem. Furthermore, Ferreau et al. used  the PAS method to solve the general convex QP problem
\begin{equation}\label{gqp}
\begin{array}{l}
\min ~\{r^{T}x+\frac{1}{2}x^{T}Q{x} |~g\leq Bx,\}
\end{array}
\end{equation}
(It is clear that the QP problem (\ref{scqp}) can be transformed into the form (\ref{gqp})) by tracking the piecewise-linear solution path of the following PQP problem
\begin{equation}\label{pqp}
\begin{array}{l}
\min ~\{r(t)^{T}x+\frac{1}{2}x^{T}Q{x} |~g(t)\leq Bx\}
\end{array}
\end{equation}
from $t=1$ to $t=0$, which has been implemented in the software package qpOASES \cite{ferreau2014qpoases}. The PQP problem in (\ref{pqp}) is constructed such that (I) when $t=0$, it becomes (\ref{gqp}), that is, $r(0)=r, g(0)=g$; (II) when $t=1$, its solution $x(1)$  as well as the corresponding multipliers $\lambda(1)\ge 0$ with $g(1)\le Bx(1)$, $\lambda(1)^T (g(1)-Bx(1))=0$ and $r(1)=-Qx(1)+B^T\lambda(1)$ are known.

At every step of PAS, it needs to solve the linear systems
\begin{equation}\label{paslinearsystem}
\left[\begin{array}{cc}Q & B_{\mathcal{A}}^T\\ B_{\mathcal{A}}&0\end{array}\right] \left[\begin{array}{c}x(t) \\\lambda_{\mathcal{A}}(t)\end{array}\right]=\left[\begin{array}{c}-r(t) \\g_{\mathcal{A}}(t)\end{array}\right],
\end{equation}
which are derived from the Karush-Kuhn-Tucker (KKT) conditions, where $\mathcal{A}$ denotes the active set. AS methods change the active set along the descent directions, while PAS changes along the parameter $t$ from $t=1$ to $t=0$. The number of steps of the PAS is generally smaller than that in the AS method. In fact, the number of  steps in the PAS  is close to the number of the different members between the starting active set and the target active set. Therefore,  the efficiency of the PAS method depends on the number of active constraints which determines the size of (\ref{paslinearsystem}), and the difference between  the starting active set and the target active set. When the starting active set and the target active set are close, PAS  needs a small number of steps to obtain the exact solution. Therefore, a good warm-wart technique for the PAS method is very important to the high performance of PAS. Compared with AS methods, the number of  steps in the PAS method is not affected by the distribution of the eigenvalues of $Q$. Specifically, if $Q$ has only a small number of large eigenvalues and if the other eigenvalues are close to zero, AS methods may be inefficient.

Because the main computation of the ALM is to solve the subproblems (\ref{alm subproblem})  (which are  special cases of (\ref{gqp}) with $B=I$), and because PAS is efficient at obtaining an exact solution of  (\ref{alm subproblem}) (when a good prediction of the optimal active is given), we combine  ALM and PAS to solve the QP problem (\ref{scqp}).  Benefiting from the framework of the ALM, the size of the KKT systems  in the combined method   is close to the number of   free variables at the solution. Therefore, the scale of the KKT systems in the combinational algorithm is smaller than that in the PM, AS and original PAS methods; this is especially significant when the solutions are sparse.

Furthermore, as mentioned in \cite{conn1991globally, dostal2003augmented}, it does not always need to obtain high-precision solutions of the subproblems;  thus, we use a first-order algorithm to approximately solve the augmented Lagrangian subproblems at the early stage of the augmented method. As $k$ increases,  the precision of the solutions of the subproblems is required to be higher. The first-order algorithms are generally efficient at obtaining   approximate solutions. However, they  needs substantially much more computation to achieve high-precision solutions. Therefore, we plan to use the PAS algorithm to obtain the exact solutions of the  augmented Lagrangian subproblems  at the mid to late stages.

However, the PAS method needs to retain the invertibility of the KKT systems in the tracking steps.  If the PAS method is applied to solve (\ref{scqp}) directly and if $Q$ is positive definite, an addition or  removal of a constraint may lead to a loss of invertibility.  In qpOASES, Ferreau et al. retain the invertibility in Eq. (\ref{paslinearsystem}) by exchanging an index of the active set and inactive set. Fortunately,  based on the framework of the ALM, the augmented Lagrangian subproblems have only bound constraints; therefore, if $Q$ is positive definite, then the Hessian matrices of $\mathcal{L}_{\beta}(x,\lambda^{k})$ are positive definite, which implies that the KKT systems in the homotopy tracking steps would always be invertible.  Thus, we do not need to exchange indices to ensure the invertibility  of the KKT systems, as is the case in qpOASES. Moreover, when $Q$ is  not positive definite, we add  proximal terms to the objective function of the augmented Lagrangian subproblems as follows
\begin{equation}\label{palm}
\begin{array}{l}
x^{k+1} =\arg\min~\{\mathcal{L}_{\beta}(x,\lambda^{k})+\frac{d_{k}}{2}\|x-x^k\|^2~|~x\geq0\}.
\end{array}
\end{equation}
Thus the Hessian matrices  of the subproblems are positive definite.

Because the PAS method is essentially a homotopy-like method for PQP, we use homotopy to denote the simplified PAS method and use AL-Hom to denote the ALM with every subproblem solved by the homotopy algorithm. Accordingly, we use PAL-Hom to denote the proximal ALM with the homotopy algorithm solving the subproblems.

Unfortunately, a simple combination of the ALM and the homotopy algorithm is unsatisfactory for QP problems (\ref{scqp}). An efficient implementation of the homotopy algorithm needs a good warm start for the homotopy algorithm as well as a fast Cholesky factor update. Moreover, although the invertibility can be ensured by the above processes, a large condition number of  (\ref{paslinearsystem}) leads to large changes in the solution, which would lead to incorrect updates of the active set. In addition, the lack of strict complementarity would also lead to incorrect updates. For these reasons, we present three important techniques: an accelerated proximal gradient method for warm starts; a sorting technique for the Cholesky factorization update, and  an $\varepsilon$-precision verification and correction technique  to correct incorrect updates of the active set.

Because the efficiency of the homotopy algorithm depends on the difference between the starting active set and the target active set, it is important to design a good warm-start technique to obtain a good estimate of the optimal active set. In the implementation of PAS, the authors have not provided methods to predict the optimal active set and just used the solution of a previous subproblem which may be not a good warm start. When PAS is directly applied to solve problems (\ref{scqp}), the performance is unsatisfactory  because for general QP problems, it is not easy to predict the active set of (\ref{scqp}). Fortunately, based on the framework of the ALM,  it is much  easier to design a warm start for PAL-Hom which iteratively solves the proximal augmented Lagrangian subproblems.

It is well known that Nesterov's accelerated proximal gradient  (APG) \cite{nesterov2005smooth, nesterov2007gradient} algorithm  is able to handle very-large-scale problems and converges at a rate $O(\frac{1}{k^2})$ which is fast for first-order algorithms. In particular, for the augmented Lagrangian subproblems, the APG is easily  implemented and has a low computational complexity  at every iteration. Moreover, the APG allows for  a rapid change in the active set at every iteration. For these reasons, we use APG to predict the optimal active set of the augmented Lagrangian subproblem. Fortunately,   a low-precision solution of (\ref{palm}) often implies a good estimate of the optimal active set; when an  approximate solution that provides a good estimate of the optimal active set is given,  the homotopy algorithm needs a small number of steps to obtain an exact solution.

Simultaneously, to improve the efficiency and  robustness of the homotopy algorithm, we present a sorting technique for the Cholesky factorization update (Section \ref{cho-fact}) that requires fewer computations than the Cholesky factorization update in qpOASES, as well as an $\varepsilon$-precision verification and a correction technique  (Section \ref{vc})  to address the incorrect updates of the active set caused by the lack of strict complementarity and the computation errors in solving the linear systems.

The outline of the remainder of this paper is as follows. Details of the homotopy algorithm are presented in Section \ref{section algorithm}. In Section \ref{applicationtolp}, we apply a proximal- point-based AL-Hom to solve linear programming (LP) problems  and prove that it converges in a finite number of iterations. Moreover, an estimate of the maximum number of  iterations and a lower bound on the descent of the linear objective are given. Finally, the numerical results for QPs and LPs from synthetic data and real-world data are presented in Section \ref{numerical}.

\section{Homotopy algorithm for the subproblems of ALM}\label{section algorithm}
In this section,  we follow Conn et al. \cite{dostal2003augmented} and Dostal et al. \cite{conn1991globally}  using the ALM to solve the QP problem (\ref{scqp}); however, we add proximal terms into the objective of the augmented Lagrangian subproblems to obtain the strict convexity of the subproblems. Moreover, because the main computation of the proximal augmented Lagrangian method is solving the augmented Lagrangian subproblems, we present a homotopy algorithm for exactly solving the   subproblems with a uniform form
\begin{equation}\label{almsub}
\begin{array}{l}
\min ~\{f^{T}z+\frac{1}{2}z^{T}H{z} | z\geq0\},
\end{array}
\end{equation}
where $H=Q+d_{k}I+\beta A^{T}A$  (if $Q$ is positive definite, $d_k=0$) is positive definite and $f=r-A^{T}\lambda^{k}-\beta A^{T}b-d_k x^k$. As mentioned in the introduction, the homotopy algorithm is a simplified PAS method, and is improved with three important techniques:  warm start, Cholesky factorization update and  $\varepsilon$-precision verification and correction.

Before presenting the homotopy algorithm, we give the optimality conditions of (\ref{almsub}) as follows, where $z^*$ is the solution of (\ref{almsub}) if and only if
\begin{eqnarray}\label{pqpkkt11}
&&Hz^*+f\geq0,\\ \label{pqpkkt12}
&&z^*\geq 0,\\\label{pqpkkt13}
&&{z^*}^T(Hz^*+f)=0,
\end{eqnarray}

\subsection{Warm start} \label{warmstart}
Because the homotopy algorithm needs a good estimate of the optimal active set of (\ref{almsub}), and because the APG  algorithm efficiently obtains an approximate solution of (\ref{almsub}), which often implies a good estimate of the optimal active set, we implement APG to approximately solve (\ref{almsub}).

Let $z^{1}=y^{0}$ which is pregiven, $l=1$, and  $\theta_{1}=1$; then, APG iterates as follows.
\begin{eqnarray}\label{apg step1}
&&y^{l}=\arg\min_{z\geq 0}\langle Hz^{l}+f,z\rangle+\frac{L}{2}\Arrowvert z-z^{l} \Arrowvert^{2}, \\	\label{apg step2}
&&\theta_{l+1}=\frac{1+\sqrt{1+4\theta_{l}^2}}{2},\\ \label{apg step3}
&&z^{l+1}=y^{l}+(\frac{\theta_{l}-1}{\theta_{l+1}})(y^{l}-y^{l-1}),
\end{eqnarray}
where $L\geq \Arrowvert H \Arrowvert$. In each iteration, $(\ref{apg step1})$ can be solved by a truncation operator
\begin{equation} \nonumber
\begin{array}{l}
y^{l}=T(z^{l}-\frac{1}{L}(Hz^{l}+f))=[z^{l}-\frac{1}{L}(Hz^{l}+f)]_{+}.
\end{array}
\end{equation}
Hence, the main computation at each iteration is a matrix-vector multiplication.

Because a low-precision solution often implies a good estimate of the optimal active set and because APG is slow at the end of the iterations, we terminate the APG algorithm when $y^{l}$ satisfies one of the following criteria.
\begin{eqnarray} \label{terminate criterion 1}
&& \mu_{\varepsilon_{1}}(y^{l}) =\mu_{\varepsilon_{1}}(y^{l-i}) ,~for~ i=1,..,S_{\max},\\\label{terminate criterion 2}
&&\frac{\Arrowvert y^{l}-y^{l-1} \Arrowvert}{\Arrowvert y^{l} \Arrowvert}<\varepsilon_{2},
\end{eqnarray}
where $\mu_{\varepsilon_{1}}(y)=\Arrowvert[y-\|y\|\varepsilon_{1}]_{+}\Arrowvert_{0}$, $S_{\max}$, $\varepsilon_{1}$ and $\varepsilon_{2}$ are some parameters that are given. Because the index $i$ is likely to be active if  $y^l_i<\eta\Arrowvert y^l \Arrowvert$, we truncate $y^l$ with $\eta\Arrowvert y^l \Arrowvert$  as follows
\begin{equation}\label{hatz}
\begin{array}{l}
\hat{z}=\left\{
\begin{array}{lcl}
y_{j}^{l},\hspace{1cm}y_{j}^{l}\geq\eta \Arrowvert y^l \Arrowvert\\
0,\hspace{1.2cm}else\\
\end{array} \right.
\end{array}
\end{equation}
where $\eta>0$. Let
\begin{eqnarray} \nonumber
 w=\left\{
\begin{array}{lcl}
-H_{j}^T\hat{z}-f_{j},      &      & {\hat{z}_{j} >  0},\\
\xi,      &      & {\hat{z}_{j} =  0},
\end{array} \right.
\end{eqnarray}
where $\xi=-\min_{j}\{H_{j}^T\hat{z}+f_{j}|\hat{z}_{j} =  0\}+\delta_1$ and $\delta_1>0$. Therefore, we have that $\hat{z}$ is the solution of
\begin{equation}\label{warm start}
\begin{array}{l}
\min ~~ \frac{1}{2}z^{T}Hz+(f+w)^{T}z \\	
\rm{s.~t.}~~~\emph{z}\geq 0.
\end{array}
\end{equation}
 from (\ref{pqpkkt11})-(\ref{pqpkkt13}).

\subsection{Homotopy tracking}

The linear homotopy between the objective function of $(\ref{almsub})$ and $(\ref{warm start})$ is
\begin{equation} \nonumber
\begin{array}{l}
h(t,z)= \frac{1}{2}z^{T}Hz+(f+tw)^{T}z~,~t\in[0,1].
\end{array}
\end{equation}

Then we can obtain the solution of $(\ref{almsub})$ by tracking the piecewise-linear solution path of the PQP problem
\begin{equation}\label{parametric qp2}
\begin{array}{l}
\min ~~ h(t,z)=\frac{1}{2}z^{T}Hz+(f+tw)^{T}z \\	
\rm{s.~t.}~~~\emph{z}\geq 0.
\end{array}
\end{equation}

Let $z(t)$, $t\in [0,1]$ be a vector function of $t$ denoting the solution path of (\ref{parametric qp2}). Suppose $z(t)$ is linear in $M$ intervals, and set $t_{0}=1, t_{M}=0$. Let $(t_{i},t_{i-1})$, $i=1,...,M$ denote the intervals, in which $z(t)$ is linear. Moreover, let $J(z(t))=\{j|H^T_j z(t)+f_{j}+tw_j=0\}$ denote the working set. Because $z(t)$ is piecewise-linear, $J(z(t))$ is constant in every interval. We use $J^{i}=\{J(z(t))|t\in(t_{i},t_{i-1})\}$ to denote the working set in the $i$-th interval, and we let $J^{i}_{c}=\{1,...,n\}\backslash J^{i}$.

\begin{proposition} \label{proposition1}
\label{th:unirep}
For any $i\in\{1,...,M\}$, there exists only one index set $S^{i}\subset \{1,...,n\}$ such that
\begin{eqnarray}\label{kkt1}
&&z_{S^{i}}(t) =-H_{S^{i}S^{i}}^{-1}(f_{S^{i}}+tw_{S^{i}})\geq 0,\\	\label{kkt2}
&&z_{S^{i}_{c}}(t)=0,\\ \label{kkt3}
&&H^{T}_{S^{i}_{c}}z(t)+f_{S^{i}_{c}}+tw_{S^{i}_{c}}>0
\end{eqnarray}
holds for any $t\in(t_{i},t_{i-1})$, where $S^{i}_{c}=\{1,...,n\}\backslash S^{i}$, $H_{S^{i}S^{i}}$ and $H_{S^{i}_{c}S^{i}}$ denote the submatrices of $H$ with appropriate rows and columns.
\end{proposition}

$Proof.$ Clearly, if $S^i$ satisfies $(\ref{kkt1})$-$(\ref{kkt3})$, then $z(t)$ is the solution of (\ref{parametric qp2}) at $t$. Moreover, because $H$ is positive definite, for any $t$, the solution of (\ref{parametric qp2}) is unique. Then, we have that $J^i$ is the unique index set, which satisfies $(\ref{kkt1})$-$(\ref{kkt3})$. \quad \endproof

The essence of the homotopy tracking steps is to calculate the solution path $z(t)$ that is unique from $t=1$ to $t=0$. This is equivalent to updating $J^i$ and $J_c^i$ from $t=1$ to $t=0$. Therefore, if a good prediction of the optimal active set is obtained, that is, $J(1)$ is close to $J(0)$, a small number of update steps is needed to change  $J(1)$ to  $J(0)$.
Benefiting from the approximate solution from APG, $\hat{z}$ is an approximate solution of $(\ref{almsub})$, and the starting active set is hopefully close to the target active set; this implies that the number of  steps of the subsequent iterations is hopefully small.

We start the homotopy tracking steps with $z(t_{0})=\hat{z}$, $J^{1}=\{j|\hat{z}_{j}>0\}$ and $J^{1}_{c}=\{1,...,n\}\backslash J^{1}$. In the homotopy tracking steps, we need to calculate $t_{i} $ and update the working set $J^{i+1}$  for $i=1,2,..M-1$.

From  Proposition \ref{proposition1}, $z(t)$ has the closed form
\begin{eqnarray} \label{closed from a}
z_{J^{i}}(t)&=&-H_{J^{i}J^{i}}^{-1}(f_{J^{i}}+tw_{J^{i}}),\\	\label{closed from b}
z_{J^{i}_{c}}(t)&=&0
\end{eqnarray}
in the $i$-th interval. We continue to decrease $t$ starting at $t_{i-1}$ until one of the following events occurs.
\begin{enumerate}[\rm{(}i)]
\item There exists $j\in J^{i}$ and $\tilde{t}<t_{i-1}$ such that $z_{j}(t)>0, t\in(\tilde{t},t_{i-1})$ and $z_{j}(\tilde{t})=0$.
\item There exists $j\in J^{i}_c$ and $\tilde{t}<t_{i-1}$ such that $H_{jJ^{i}}z_{J^{i}}(\tilde{t})+(f_{j}+\tilde{t}w_{j})=0$.
\end{enumerate}
When $($i$)$ or $($ii$)$ occurs, we need to calculate the value of $\tilde{t}$, and $J^{i}$ and $J^{i}_{c}$ need to exchange indices at $\tilde{t}$.

According to $($i$)$ and $($ii$)$, define
\begin{eqnarray}\nonumber
\hat{j}= \rm arg\max_{\emph j} &&\{\frac{u^{i}_{j}}{v^{i}_{j}}< t_{i-1}|j\in J^{i} ~and~ v^{i}_{j}<0\},\\  \nonumber
\tilde{j}=\rm arg\max_{\emph j} && \{\frac{\psi^{i}_{j}}{\phi^{i}_{j}}< t_{i-1}|j\in J^{i}_{c} ~and~ \phi^{i}_{j}<0\},
\end{eqnarray}
where $u^{i}= -H_{J^{i}J^{i}}^{-1}f_{J^{i}}$, $v^{i}= H_{J^{i}J^{i}}^{-1}w_{J^{i}}$, $\psi^{i}= H_{J^{i}_{c}J^{i}}u^{i}+f_{J^{i}_{c}}$ and $\phi^{i}= H_{J^{i}_{c}J^{i}}v^{i}-w_{J^{i}_{c}}$.

If $\hat{j}$ is empty, set $\frac{u^{i}_{\hat{j}}}{v^{i}_{\hat{j}}}=-\infty$, which is the same as $\tilde{j}$. Now, we discuss the update strategy of $J^{i}$ and $J_c^{i}$ as follows.

\noindent  $\textbf{Case 1:}$  $\frac{u^{i}_{\hat{j}}}{v^{i}_{\hat{j}}}>\frac{\phi^{i}_{\tilde{j}}}{\psi^{i}_{\tilde{j}}}$ and $\frac{u^{i}_{\hat{j}}}{v^{i}_{\hat{j}}}>0$.

\noindent Thus, $($i$)$ occurs first; then, we obtain $t_{i}=\tilde{t}=\frac{u^{i}_{\hat{j}}}{v^{i}_{\hat{j}}}$, $J^{i+1}=J^{i}\backslash \hat{j}$ and $J^{i+1}_{c}=J^{i}_{c}\cup \hat{j}$. Thus, $z(t)$ has the following closed form
\begin{eqnarray}\label{closed from1}
z_{J^{i+1}}(t) &=&-H_{J^{i+1}J^{i+1}}^{-1}(f_{J^{i+1}}+tw_{J^{i+1}}),\\	\label{closed from2}
z_{J^{i+1}_{c}}(t)&=&0
\end{eqnarray}
in the interval $(t_i,t_{i+1})$.

\noindent  $\textbf{Case 2:}$  $\frac{u^{i}_{\hat{j}}}{v^{i}_{\hat{j}}}<\frac{\phi^{i}_{\tilde{j}}}{\psi^{i}_{\tilde{j}}}$ and $\frac{\phi^{i}_{\tilde{j}}}{\psi^{i}_{\tilde{j}}}>0$.

\noindent Thus, $($ii$)$ occurs first; then, $t_{i}=\tilde{t}=\frac{\phi^{i}_{\tilde{j}}}{\psi^{i}_{\tilde{j}}}$, $J^{i+1}=J^{i}\cup \tilde{j}$ and $J^{i+1}_{c}=J^{i}_{c}\backslash \tilde{j}$.

\noindent  $\textbf{Case 3:}$  $\frac{u^{i}_{\hat{j}}}{v^{i}_{\hat{j}}}\leq 0$, $\frac{\phi^{i}_{\tilde{j}}}{\psi^{i}_{\tilde{j}}}\leq 0$.

\noindent In this case, the algorithm will terminate and we obtain

\begin{equation}\label{endsolution}
\begin{array}{l}
z_{J^{i}}(0) =-H_{J^{i}J^{i}}^{-1}f_{J^{i}},\\ 	
z_{J^{i}_{c}}(0)=0.
\end{array}
\end{equation}

Note that  $w$ is constructed such that $\hat{z}$ satisfies the strict complementarity conditions  at $t=1$. However, in the homotopy tracking steps, there may exist an interval $(t_{i}, t_{i-1})$ such that for some $j_1$
\begin{eqnarray} \label{complementary}
H_{j_1}^T z(t)+f_{j_1}+tw_{j_1}=0~{\rm and}~z_{j_1}(t)=0, t\in (t_{i}, t_{i-1}).
\end{eqnarray}
Because the above update strategy does not consider these indices, we need to check whether (\ref{complementary}) still holds with $J^i$ and $J^i_c$ exchanging indices as above in the $(i+1)$-th interval. Specifically, we simply need to check the value of $v^{i+1}_{j_1}$. If $v^{i+1}_{j_1}>0$, then the strictly complementarity conditions hold at the ${j_1}$-th component; if  $v^{i+1}_{j_1}=0$, then the strictly complementarity conditions do not hold, and if $v^{i+1}_{j_1}<0$, then we add ${j_1}$ to $J^{i+1}_c$.

By tracking the solution path of $(\ref{parametric qp2})$ as above, we obtain $\bar{z}=z(0)$, which is the solution of (\ref{almsub}).

Clearly, the complexity of the homotopy algorithm depends on the number of the steps and the size of $J^i$. Specifically, at every step, we need to solve two symmetric positive-definite linear systems of equations
\begin{eqnarray} \label{linear equation}
H_{J^{i}J^{i}}{u^{i}}= f_{J^{i}}~~and~~H_{J^{i}J^{i}}{v^{i}}= w_{J^{i}}
\end{eqnarray}
and perform one matrix-vector multiplication
\begin{eqnarray} \label{matrix vector product}
[\psi^{i},\phi^{i}]= H_{J^{i}_{c}J^{i}}[{u^{i}},{v^{i}}]+[f_{J^{i}_{c}},-w_{J^{i}_{c}}].
\end{eqnarray}
We simply need to solve one equation in (\ref{linear equation}) for
\begin{eqnarray} \nonumber
u^{i}+t_{i-1}v^{i}=x_{J^{i}}(t_{i-1}).
\end{eqnarray}
Thus, when $|J^{i}|$ is small, the homotopy algorithm has a low computational complexity  at each step. In addition, benefiting from the approximate solution from APG, the number of the steps is hopefully   small.

Unlike the original PAS method  for QP problem (\ref{scqp}), the PAL-Hom algorithm would always ensure the strict convexity for both adding and removing an index; therefore, we do not need to check the invertibility after exchanging an index. Because $H_{J^{i}J^{i}}$ is positive definite, we apply the Cholesky factorization method for (\ref{linear equation}). Moreover, because $J^{i}$ changes one member  every time and because the exchanged index is more likely to be the index whose corresponding value is close to zero, we present a  sorting technique for the Cholesky factorization update  different from that in qpOASES  \cite{ferreau2006online,ferreau2014qpoases}.

\subsection{Update the Cholesky factorization} \label{cho-fact} Note that the index $j_{1}$ is more likely to be active than $j_{2}$ if $\hat{z}_{j_{1}}<\hat{z}_{j_{2}}$; thus $j_{1}$ is more likely to be removed from $J^{i}$ than $j_{2}$  in the homotopy tracking steps. For this reason, at the start of the homotopy tracking steps, we sort $J(\hat{z})$ by the value of $\hat{z}_{j}, j\in J(\hat{z})$, that is,
$$\hat{z}_{[J(\hat{z})]_{s}}\geq\hat{z}_{[J(\hat{z})]_{s+1}},$$
where $[J(\hat{z})]_{s}$ denotes the $s$-th member of $J(\hat{z})$. With this sorting technique, the indices corresponding to the smaller $\hat{z}_{j}$  would be sorted at the end of $J(\hat{z})$; thus, the indices removed from $J^{i}$ would be distributed at the end of $J^{i}$. Moreover, when an index is added to $J^{i}$, we put it at the end of $J^{i}$.

 Assume that $J^{i}$ is known and that $H_{J^{i}J^{i}}$ has the Cholesky factorization
$$R^{T}R=H_{J^{i}J^{i}}.$$
Then we update the  Cholesky factorization as follows.

\hangafter=1
\setlength{\hangindent}{1.02cm}
$\rhd$ Add an index $\tilde{j}$ to $J^{i}$; then,
$$H_{J^{i+1}J^{i+1}}=\left[\begin{array}{cc}H_{J^{i}J^{i}}&H_{J^{i}\tilde{j}}\\H_{\tilde{j}J^{i}}&H_{\tilde{j}\tilde{j}}\end{array}\right].$$
Let $H_{J^{i+1}J^{i+1}}=\tilde{R}^{T}\tilde{R}$ be the Cholesky factorization; then,
$$\tilde{R}=\left[\begin{array}{cc}R&\tilde{r}\\0&\sqrt{H_{\tilde{j}\tilde{j}}-\tilde{r}^{T}\tilde{r}}\end{array}\right],$$
where $R^{T}\tilde{r}=H_{J^{i}\tilde{j}}$. This update requires only $\frac{1}{2}\Gamma_{i}^2$ flops, where $\Gamma_{i}=|J^{i}|$.

\hangafter=1
\setlength{\hangindent}{1.02cm}
$\rhd$ Remove an index $\hat{j}$ from $J^{i}$, then
$$H_{J^{i+1}J^{i+1}}=\left[\begin{array}{cc}H_{J_{1}^{i}J_{1}^{i}}&H_{J_{1}^{i}J_{2}^{i}}\\H_{J_{2}^{i}J_{1}^{i}}&H_{J_{2}^{i}J_{2}^{i}}\end{array}\right],$$
where $J^{i}=[J_{1}^{i},\hat{j},J_{2}^{i}]$. Assume that $H_{J^{i+1}J^{i+1}}=\hat{R}^{T}\hat{R}$ is the Cholesky factorization; then, we have
$$\hat{R}=\left[\begin{array}{cc}R_{I^{i}_{1}I^{i}_{1}}&R_{I^{i}_{1}I^{i}_{2}}\\0&\bar{R}\end{array}\right],$$
where $I^{i}_{1}=\{1,...,|J^{i}_{1}|\}$, $I^{i}_{2}=\{|J^{i}_{1}|+2,...,|J^{i}|\}$ and $\bar{R}^{T}\bar{R}=H_{J_{2}^{i}J_{2}^{i}}-R_{I^{i}_{1}I^{i}_{2}}^{T}R^{i}_{I^{i}_{1}I^{i}_{2}}$. This case therefore requires $\frac{2}{3}|J^{i}_{2}|^3$ flops.

In conclusion
\begin{equation}\label{equation4.20_1a}
\begin{array}{l}
\left\{
\begin{array}{lcl}
\frac{1}{2}\Gamma_{i}^2,\hspace{3.2cm}add;\\
\frac{2}{3}|J^{i}_{2}|^3+(\Gamma_{i}-|J^{i}_{2}|)|J^{i}_{2}|^2,\hspace{0.15cm}remove;\\
\end{array} \right.
\end{array}
\end{equation}
flops are required to update the Cholesky factorization at each step, where $(\Gamma_{i}-|J^{i}_{2}|)|J^{i}_{2}|^2$ is the matrix multiplication $R_{I^{i}_{1}I^{i}_{2}}^{T}R_{I^{i}_{1}I^{i}_{2}}$, while the Cholesky factorization update technique \cite{ferreau2006online} of the PAS method in qpOASES would require
\begin{equation}\label{equation4.20_1b}
\begin{array}{l}
\left\{
\begin{array}{lcl}
5\Gamma_{i}^2,\hspace{1.9cm}add;\\
\frac{5}{2}\Gamma_{i}^2,\hspace{1.9cm}remove;\\
\end{array} \right.
\end{array}
\end{equation}
flops at each step. Our update strategy requires fewer computations when adding an index than that in qpOASES. Moreover, benefiting from the sorting technique, $|J^{i}_{2}|\ll \Gamma_{i}$; therefore, the removing update is a low-cost technique.

\subsection{ $\varepsilon$-precision verification and correction}\label{vc} From the homotopy tracking steps, we have
\begin{eqnarray} \label{close_formJi1}
&&z_{J^{i}}(t_i)=-H_{J^{i}J^{i}}^{-1}(f_{J^{i}}+tw_{J^{i}}),\\	\label{close_formJi2}
&&z_{J^{i}_{c}}(t_i)=0.
\end{eqnarray}
However, due to the errors from the solving of the linear systems which may have a large condition number, the update of $J^i$ and $J^i_c$ may not be correct. Moreover, the lack of strict complementarity may also lead to an incorrect update of $J^i$ and $J^i_c$; therefore, we need to verify that  $z(t_i)$ satisfies  the optimality conditions.
\begin{eqnarray} \label{kkt at ti}
&&z_{J^{i}}(t_i)\geq0, \\	\label{kkt at ti1}
&&H_{J_{c}^{i}J^{i}}z_{J^{i}}(t_i)+f_{J^{i}}+tw_{J^{i}}\geq0,
\end{eqnarray}
In practice, it is not necessary and may  be difficult to ensure that (\ref{kkt at ti})-(\ref{kkt at ti1}) strictly hold, especially when the strict complementarity conditions are weak; therefore, we relax (\ref{kkt at ti})-(\ref{kkt at ti1}) by a small $\varepsilon$ as follows.
\begin{eqnarray} \label{kkt at ti2}
&&z_{J^{i}}(t_i)\geq -\varepsilon, \\	\label{kkt at ti3}
&&H_{J_{c}^{i}J^{i}}z_{J^{i}}(t_i)+f_{J^{i}}+tw_{J^{i}}\geq -\varepsilon,
\end{eqnarray}
If (\ref{kkt at ti2})-(\ref{kkt at ti3}) hold, the homotopy algorithm goes to the next step; otherwise, we correct $J^{i}$ and $J^{i}_{c}$ as follows.

\textbf{Step 1}: If there exists $j\in J^{i}$ such that $z_j(t_i)<-\varepsilon$, then let
$$\bar{j}=\arg\min_{j\in J^i} ~\{z_j(t_i)\}$$
and $J^{i}=J^{i}\cup \bar{j}, J_{c}^{i}=J_{c}^{i}\backslash \bar{j}$, refresh $z(t_i)$ as in (\ref{close_formJi1})-(\ref{close_formJi2}) and go to \textbf{Step 1}; otherwise go to \textbf{Step 2}.

\textbf{Step 2}: If there exists $j\subset J_{c}^{i}$ such that $H^T_{j}z(t_i)+f_{j}+t_{i}w_{j}\leq -\varepsilon$, then let
$$\bar{j}=\arg\min_{j\in J^i_c} ~\{H^T_{j}z(t_i)+f_{j}+t_{i}w_{j}\}$$
and $J_{c}^{i}=J_{c}^{i}\backslash \bar{j}, J^{i}=J^{i}\cup \bar{j}$, refresh $z(t_i)$ as in (\ref{close_formJi1})-(\ref{close_formJi2}) and go to \textbf{Step 1}; otherwise, terminate the correction steps.

The correction steps ensure that the solution $x(t)$ satisfies the optimality conditions with $\varepsilon$-precision and guarantee the stability of the homotopy tracking algorithm.

Finally, as mentioned in the introduction, in many cases, it is not necessary to obtain the exact solutions of the first few augmented Lagrangian subproblems; therefore, for these subproblems, we directly go to the next iteration after the approximate solution is obtained by APG. For the other subproblems, we use the homotopy algorithm to obtain exact solutions. The framework of PAL-Hom for convex QP is given as Algorithm \ref{ALM algorithm}.
\begin{algorithm}[!h]
    \caption{PAL-Hom algorithm for  QP}
    \begin{algorithmic}[]
    \REQUIRE ~~\\                          %Parameters Input
     $k=0$, $x^{0}$, $\lambda^{0}$, $\beta$, $tol$, $\varepsilon_c>0$
    \ENSURE ~~\\
      $x^{k+1}$;

      \WHILE{$\|Ax^{k}-b\|>tol$ or $\|x^{k}-x^{k-1}\|>tol$ }
       \STATE Approximately solve $(\ref{palm})$ with APG algorithm as $(\ref{apg step1})$, $(\ref{apg step2})$ and $(\ref{apg step3})$ until $(\ref{terminate criterion 1})$ or $(\ref{terminate criterion 2})$ is satisfied.
       \IF{$\|Ax^{k}-b\|<\varepsilon_{c}$}
       \STATE Track the solution path of $(\ref{parametric qp2})$ from $t=1$ to $t=0$ and set $x^{k+1}$ to $z(0)$ in (\ref{endsolution}).
       \ELSE
       \STATE Set $x^{k+1}$ to $\hat{z}$ in (\ref{hatz}).
       \ENDIF

       \STATE $\lambda^{k+1} =\lambda^{k}-\beta(Ax^{k+1}-b)$;
       \STATE $k=k+1$;
      \ENDWHILE

    \end{algorithmic} \label{ALM algorithm}

\end{algorithm}

\section{An application to LP}\label{applicationtolp}
Because LP
\begin{equation}\label{lp}
\begin{array}{l}
\min ~~ c^{T}x ~~~~\\ 	
\rm{s.t.}~~~\emph{Ax}=\emph{b},\\
~~~~~~~\emph{x}\geq 0,
\end{array}
\end{equation}
is a special case of QP with $Q=\bf{0}$ and $r=c\in \mathbb{R}^n$,  PAL-Hom can be applied to solve LP problems. Moreover, for LP problems, Wright \cite{wright1990implementing} showed that PAL-Hom converges in a finite number of steps if the subproblems are exactly solved  for all $k$ sufficiently large and if the strict complementarity conditions hold at the solution.

On the other hand,  Mangasarian  \cite{mangasarian1981iterative, mangasarian1979nonlinear} transformed the LP problem into a weakly strictly convex QP problem
\begin{equation}\label{regularized lp}
\begin{array}{l}
\min ~~ c^{T}x+\frac{\varepsilon}{2}x^{T}{x} ~~~~\\ 	
\rm{s.t.}~~~\emph{Ax}=\emph{b},\\
~~~~~~~x\geq 0
\end{array}
\end{equation}
by adding a small regularization term to the objective. Moreover, Mangasarian proved that (\ref{regularized lp}) obtains a solution of (\ref{lp}) if $\varepsilon$ is smaller than some $\bar{\varepsilon}>0$.  However, it is difficult to derive a realistic priori estimate of $\bar{\varepsilon}$, and for certain practical problems, $\bar{\varepsilon}$ would be very small. If we apply AL-Hom to solve (\ref{regularized lp}),  a small $\bar{\varepsilon}$ would lead to a large condition number of the KKT systems in the homotopy tracking steps, which is adverse to the robustness of the homotopy algorithm.

Motivated by Mangasarian \cite{mangasarian1981iterative}, we used proximal point methods to solve LP problems, that is, for a given $x^0\in \mathbb{R}^n$, iteratively solve the strictly convex subproblems
\begin{eqnarray}\label{PPA}
x^{\sigma+1}=\min_{x\in \Omega}~ c^{T}x+\frac{1}{2\alpha_{\sigma}}\|x-x^\sigma\|^2, \sigma=0,1,2,....
\end{eqnarray}
Moreover, every subproblem is solved by AL-Hom. We use  PP-AL-Hom to denote the  above process for LP.

Under the assumption that $(\ref{lp})$ has at least one finite solution, we prove that,  if $\alpha^{\sigma}>\alpha$ for some $\alpha>0$, then iterations (\ref{PPA}) converge in a finite number of steps. Simultaneously, we give a positive lower bound of $c^{T}x^{\sigma}-c^{T}x^{\sigma+1}$ and an estimate of the maximum number of the iterations  (\ref{PPA}). Since $\alpha$ can be arbitrary, the condition number of the KKT systems in the homotopy tracking steps can be controlled.

It is  clear that (\ref{regularized lp}) is a  special case of (\ref{PPA}) with $\alpha^{0}=\frac{1}{\varepsilon}$ and $x^0=0$. Therefore, we have that $(\ref{PPA})$ converges in one step, if $\alpha^{0}\geq\frac{1}{\bar{\varepsilon}}$ and $x^0=0$. Moreover, in contrast to PAL-Hom, the finite-step termination of PP-AL-Hom does not require the strictly complementarity conditions at the solution.

\begin{figure}[!htb]
  \centering
  \vspace{-3.5cm}
  \includegraphics[width=9.5cm]{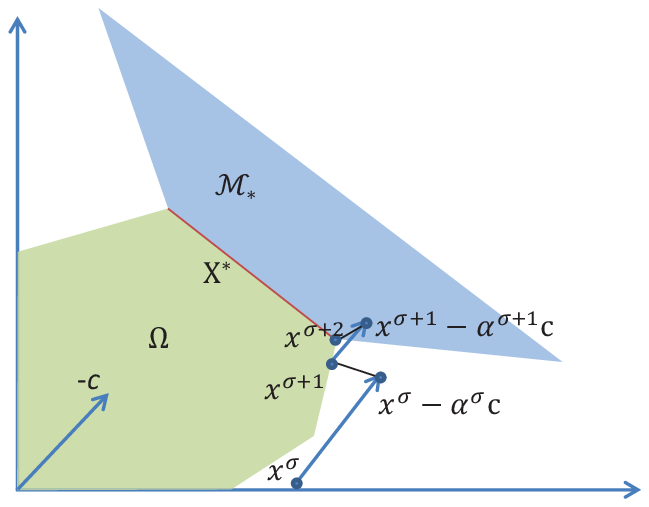}\\
  \vspace{-5cm}
  \caption{Projected gradient method for LP}
\label{figure1}
\end{figure}

Let $\Omega=\{x|Ax=b, x\geq0\}$ and $X^{*}$ denote the solution set of $(\ref{lp})$. Define
\begin{eqnarray}\nonumber
\mathcal{M}_{*}=\bigcup_{x^{*}\in X^{*}} (x^{*}+N_{\Omega}(x^*)),
\end{eqnarray}
where $N_{\Omega}(x^{*})$ is the normal cone of $\Omega$ at $x^*$. Clearly, (\ref{PPA}) is equivalent to
\begin{eqnarray}\label{pg subproblem}
x^{\sigma+1}=P_{\Omega}(x^{\sigma}-\alpha^{\sigma} c),
\end{eqnarray}
where  $P_{\Omega}(y)=\arg\min_{x\in \Omega} \frac{1}{2}\|x-y\|^2$ is the projection operator onto $\Omega$. Therefore, (\ref{PPA}) is equivalent to the projection procedures in Figure \ref{figure1}. It is clear that if $x^{\sigma}-\alpha^{\sigma}$ is local in $\mathcal{M}_{*}$,  then $x^{\sigma+1}$ is the solution of (\ref{lp}).

%%%%%%%%%%%%%%%%%%%%%%%%%
%%%%%%%%%%%%%%%%%

\begin{theorem}
\label{th:gradientint}
 $-c\in$ \rm{int} $\mathcal{M}^{\infty}_{*}$, \emph{where} $\mathcal{M}^{\infty}_{*}$ \emph{is the asymptotic cone of}  $\mathcal{M}_{*}$ and $\rm{int}~\mathcal{M}^{\infty}_{*}$ denotes the interior of  $\mathcal{M}^{\infty}_{*}$.
\end{theorem}

$Proof.$ We know from the optimality conditions that
\begin{equation}\nonumber
\begin{array}{l}
-c\in N_{\Omega}(x^{*}), \forall x^{*}\in X^{*}.
\end{array}
\end{equation}
Then we have $x^{*}+t(-c)\in \mathcal{M}_{*}$, $\forall t\geq 0$, which implies
\begin{equation}\label{norm cone condition}
\begin{array}{l}
-c\in \mathcal{M}^{\infty}_{*}.
\end{array}
\end{equation}

Define
\begin{equation}\nonumber
\begin{array}{l}
\mathcal{T}(x^{*},d)=\{t|x^{*}+td\in X^{*}\},
\end{array}
\end{equation}
where $d\in R^{n}$ and satisfies $d^{T}c=0$. Because $X^{*}$ is a closed convex set, $\mathcal{T}(x^{*},d)$ is a closed interval. Next, we prove that
\begin{equation}\label{norm cone condition in subspace}
\begin{array}{l}
-c \in \rm ri~(\rm{span}\{-\emph{c},\emph{d}\}\cap \mathcal{M}^{\infty}_{*}),
\end{array}
\end{equation}
where ri $S$ denotes the relative interior of $S$.

If $\mathcal{T}(x^{*},d)=(-\infty,+\infty)$, clearly, $d$, $-d\in \mathcal{M}^{\infty}_{*}$. Therefore, (\ref{norm cone condition in subspace}) is obvious by (\ref{norm cone condition}).

If $\mathcal{T}(x^{*},d)=(-\infty,t_{\max}]$ and $t_{\max}<\infty$. Similar to above, $-d\in \mathcal{M}^{\infty}_{*}$. Moreover, there exists $u\in N_{\Omega}(x^{*}+t_{\max}d)\cap \rm{span}\{-\emph{c},\emph{d}\}$ that satisfies $\langle u,-d \rangle <0$; therefore, (\ref{norm cone condition in subspace}) holds for $d^{T}c=0$.  Moreover, if there exists no such $u$, we can find $t^{'}_{\max}>t_{\max}$ such that $x^{*}+t^{'}_{\max}d \in X^{*}$ because $\Omega$ is a convex polyhedron, which contradicts the definition of $\mathcal{T}(x^{*},d)$.

If $\mathcal{T}(x^{*},d)=[t_{\min},+\infty)$ and $t_{\min}>-\infty$. let $\tilde{d}=-d$; then, $\mathcal{T}(x^{*},\tilde{d})=(-\infty,-t_{\min}]$. Therefore, we have $-c \in \rm ri~(\rm{span}\{-\emph{c},\tilde{\emph{d}}\}~\cap~ \mathcal{M}^{\infty}_{*})=\rm ri~(\rm{span}\{-\emph{c},\emph{d}\}~\cap~ \mathcal{M}^{\infty}_{*})$ from the situation above.

If $\mathcal{T}(x^{*},d)=[t_{\min},t_{\max}]$ and $t_{\min}>-\infty$, $t_{\max}<\infty$, then, similar to above, there exists $u_{1}\in N_{\Omega}(x^{*}+t_{\min}d)~\cap~ \rm{span}\{-\emph{c},\emph{d}\}$ and $u_{2}\in N_{\Omega}(x^{*}+t_{\max}d)~\cap~ \rm{span}\{-\emph{c},\emph{d}\}$ such that
\begin{eqnarray}\label{equation2.6}
\langle u_{1},d \rangle <0~\rm and~\langle \emph u_{2},\emph d \rangle >0. \nonumber
\end{eqnarray}
Thus, (\ref{norm cone condition in subspace}) holds for $d^{T}c=0$.

Because $d$ is arbitrary in the space $\{s|s^{T}c=0\}$, we have $-c \in \rm int\hspace{0.03cm}(\mathcal{M}^{\infty}_{*})$ from (\ref{norm cone condition in subspace}). ~~~ \endproof

\begin{theorem}
\label{th:lptoqp}
For any $x^{0}\in \mathbb{R}^{n}$, assume that the sequence $\{x^{\sigma}\}$ is obtained by $(\ref{pg subproblem})$; then,
\begin{enumerate}[\rm{(}i)]
\item There exists an $\bar{\alpha}$ such that if $ \alpha^0\geq\bar{\alpha}$, then $x^{1}\in X^{*}$.
\item There exists $0<\theta_{\min}\leq\frac{\pi}{2}$ such that
\begin{equation}\nonumber
\begin{array}{l}
c^{T}x^{\sigma}-c^{T}x^{\sigma+1}\geq \alpha^{\sigma} (1-\cos\hspace{0.05cm}\theta_{\min})\Arrowvert c \Arrowvert^{2},
\end{array}
\end{equation}
if $x^{\sigma+1}\notin X^{*}$. Moreover, $\theta_{\min}=\arccos\left(\frac{\Arrowvert P_{\rm bd(\mathcal{M}^{\infty}_{*})}(-c) \Arrowvert}{\Arrowvert c \Arrowvert}\right)$, where $\rm bd(\hspace{0.03cm}\mathcal{M}^{\infty}_{*})$ denotes the boundary of $\mathcal{M}^{\infty}_{*}$, $P_{\rm bd(\mathcal{M}^{\infty}_{*})}(\cdot)$ denotes the projection onto {\rm bd}$(\mathcal{M}^{\infty}_{*})$.
\item For any $\alpha>0$, $p\in \{-1 \cup  \mathcal{N}^{+}\}$, if $\alpha^{\sigma}\geq\alpha$, for $\sigma=p,p+1,...$, then there exists $\Gamma=\left[\frac{c^{T}x^{p+1}-c^{T}x^{*}}{\alpha (1-\cos\hspace{0.05cm}\theta_{\min})\Arrowvert c \Arrowvert^{2}}+p+2\right]_{+}$, such that
$x^{\Gamma}\in X^{*}$.
\end{enumerate}
\end{theorem}

$Proof.$ We prove each of the three claims in turn.

 (i) Define $B(r)=\{x|\Arrowvert x\Arrowvert\leq r\}$. Because $-c\in$ int $ \mathcal{M}^{\infty}_{*}$, there exists $\varepsilon>0$ such that
$$-c+B(\varepsilon)\subset \mathcal{M}^{\infty}_{*}.$$
Then for any $x\in \mathcal{M}_{*}$ and $\alpha>0$, we arrive at
$$x+\alpha B(\varepsilon)-\alpha c=x+\alpha (-c+B(\varepsilon))\subset \mathcal{M}_{*}.$$
Let $\bar{\alpha}=\frac{\Arrowvert x^{0}-x\Arrowvert}{\varepsilon}$; hence,
\begin{eqnarray}\label{equation2.9} \nonumber
x^{0}-\bar{\alpha}c&=&x+x^{0}-x-\bar{\alpha}c \\ \nonumber
                   &=&x+\frac{\Arrowvert x^{0}-x\Arrowvert}{\varepsilon} \cdot \frac{\varepsilon(x^{0}-x)}{\Arrowvert x^{0}-x\Arrowvert}-\bar{\alpha}c \\ \nonumber
                   &\subset& x+\bar{\alpha}B(\varepsilon)-\bar{\alpha}c\\\nonumber
                   &\subset& \mathcal{M}_{*}.
\end{eqnarray}
Moreover, for any $\alpha\geq \bar{\alpha}$
\begin{eqnarray}\label{equation2.10} \nonumber
x^{0}-\alpha c&=&x^{0}-\bar{\alpha}c -(\bar{\alpha}-\alpha)c \\ \nonumber
              &\subset& -(\bar{\alpha}-\alpha)c+ \mathcal{M}_{*}\\ \nonumber
              &\subset& \mathcal{M}_{*}.
\end{eqnarray}
Thus, from the definition of $\mathcal{M}_{*}$, we have $x^{1}=P_{\Omega}(x^{0}-\alpha^0 c)\in X^{*}$,  $\alpha^0\geq\bar{\alpha}$ .

(ii) If $x^{\sigma+1}=x^{\sigma}-\alpha^{\sigma} c$, then
$$c^{T}x^{\sigma}-c^{T}x^{\sigma+1}= \alpha^{\sigma} \Arrowvert c \Arrowvert^{2}.$$

If $x^{\sigma+1}\neq x^{\sigma}-\alpha^{\sigma} c$, then
$$\langle x^{\sigma}-x^{\sigma+1},x^{\sigma}-\alpha^{\sigma} c-x^{\sigma+1}\rangle\geq \frac{\pi}{2}$$
holds by (\ref{pg subproblem}), where $\langle s_{1},s_{2}\rangle$ denotes the angle between $s_{1}$ and $s_{2}$. Let $\theta$ denote the angle between $x^{\sigma}-\alpha^{\sigma} c -x^{\sigma+1}$ and $-c$. Because $x^{\sigma+1}\notin X^{*}$, we obtain $x^{\sigma}-\alpha^{\sigma} c -x^{\sigma+1}\notin \rm int \mathcal{M}^{\infty}_{*}$ from the convexity of $\Omega$. Thus, we have $\theta\geq \theta_{\min}=\arccos(\frac{\Arrowvert P_{\rm [bd\mathcal{M}^{\infty}_{*}]}(-c) \Arrowvert}{\Arrowvert c \Arrowvert})$ from $-c\in \rm int  \hspace{0.03cm}\mathcal{M}^{\infty}_{*}$. So
\begin{eqnarray} \nonumber
c^{T}x^{\sigma}-c^{T}x^{\sigma+1}&=&c^{T}x^{\sigma}-c^{T}P_{[x^{\sigma},x^{\sigma}-\alpha^{\sigma} c]}(x^{\sigma+1})\\\label{value decrease bound}
                       &=&c^{T}x^{\sigma}-c^{T}(x^{\sigma}-\alpha^{\sigma}(1-\cos\hspace{0.03cm}\theta)c) \\\nonumber
                       &\geq& \alpha^{\sigma} (1-\cos\hspace{0.05cm}\theta_{\min})\Arrowvert c \Arrowvert^{2},
\end{eqnarray}
where $[x^{\sigma},x^{\sigma}-\alpha^{\sigma} c]$ denotes a segment whose endpoints are $x^{\sigma}$ and $x^{\sigma}-\alpha^{\sigma} c$.  Note that because $-c\in$ \rm{int} $\mathcal{M}^{\infty}_{*}$, we have $1-\cos\hspace{0.05cm}\theta_{\min}>0$.

(iii) For any $\kappa\geq 1$, if $x^{\sigma}\notin X^{*},~\sigma=1,2,..,\kappa$, then we have from (\ref{value decrease bound}) that
\begin{eqnarray}\nonumber
c^{T}x^{p}-c^{T}x^{*}&\geq& c^{T}x^{p}-c^{T}x^{\kappa}\\\label{decrease}
                       &=&c^{T}x^{p}-c^{T}x^{p+1}+\sum_{\sigma=p+1}^{\kappa-1}(c^{T}x^{\sigma}-c^{T}x^{\sigma+1})\\\nonumber
                       &\geq&c^{T}x^{p}-c^{T}x^{p+1}+(\kappa-p-1)\alpha (1-\cos\hspace{0.05cm}\theta_{\min})\Arrowvert c \Arrowvert^{2},
\end{eqnarray}
  which implies
\begin{eqnarray}\nonumber
\kappa\leq\frac{c^{T}x^{p+1}-c^{T}x^{*}}{\alpha (1-\cos\hspace{0.05cm}\theta_{\min})\Arrowvert c \Arrowvert^{2}}+p+1.
\end{eqnarray}
Then we have (iii) when
\begin{eqnarray}\nonumber
\Gamma=\left[\frac{c^{T}x^{p+1}-c^{T}x^{*}}{\alpha (1-\cos\hspace{0.05cm}\theta_{\min})\Arrowvert c \Arrowvert^{2}}+p+2\right]_{+}.
\end{eqnarray}
 ~~~\endproof

Moreover, from (\ref{decrease}), we have
\begin{eqnarray}\nonumber
c^{T}x^{p+1}-c^{T}x^{*}&\geq& \sum_{\sigma=p+1}^{\kappa-1}\alpha^{\sigma} (1-\cos\hspace{0.05cm}\theta_{\min})\Arrowvert c \Arrowvert^{2}\\\nonumber
\end{eqnarray}
if $x^{\kappa}\notin X^{*}$, which implies
\begin{eqnarray}\nonumber
\sum_{\sigma=p+1}^{\kappa-1}\alpha^{\sigma}\leq \frac{c^{T}x^{p+1}-c^{T}x^{*}}{(1-\cos\hspace{0.05cm}\theta_{\min})\Arrowvert c \Arrowvert^{2}}.\\\nonumber
\end{eqnarray}
Then we obtain $x^{\kappa}\in X^{*}$ so long as
\begin{eqnarray}\nonumber
\sum_{\sigma=p+1}^{\kappa-1}\alpha^{\sigma}\geq \frac{c^{T}x^{p+1}-c^{T}x^{*}}{(1-\cos\hspace{0.05cm}\theta_{\min})\Arrowvert c \Arrowvert^{2}}.\\\nonumber
\end{eqnarray}

By Theorem \ref{th:lptoqp}, we have that if LP problem $(\ref{lp})$ has one finite solution, it can be transformed into a finite number of strictly convex QP problems with projection form similar to $(\ref{pg subproblem})$, which is equivalent to (\ref{PPA}). We solve every projection problem by using AL-Hom.

\section{Numerical results}\label{numerical}

In this section, we demonstrate the performance of our algorithms. The numerical experiments were performed on the MATLAB 8.1 programming platform (R2013a) running on a machine with the a Windows 7 operating system, an Intel(R) Core(TM)i7 6700 3.40GHz processor and 32 GB of memory. The QP-solvers and LP-solvers in the other software packages were called by the MATLAB interface.

We tested PAL-Hom for solving randomly generated QPs  and QPs from the CUTEr test set\cite{bongartz1995cute}. We also used PAL-Hom to solve the discrete contact problems of elasticity and QPs from SVMs \cite{sra2012optimization} that were applied to speech recognition and handwritten digit recognition. Finally, we used PP-AL-Hom to solve randomly generated LPs and LPs from the Netlib test set \cite{gay1985electronic}.

\subsection{Experiments on QPs from synthetic data and CUTEr test set}
\textbf{$\bullet$ Randomly generated QPs.} In this part, we randomly generated dense and sparse standard QPs $(\ref{scqp})$ with MATLAB codes as follows.

~~~~ A=sprandn($m,n,d_A$); B=sprandn($q,n,d_{B}$); Q=B$'*$B;

~~~~~r=$-$B$'*$randn($m$,1); b=$10*$randn($m$,1),

\noindent where $d_A, d_{B}$ denote the density of $A$ and $B$ which are pregiven, $d_Q$  denotes the density of $Q$, and ``randn" denotes normally random distribution  function.

\begin{table}[!htb]
\vspace{-0.4cm}
\caption{\scriptsize{Randomly generated QPs}}
 \tiny
\begin{tabular}{p{1.7cm}p{1.2cm}p{1.2cm}p{1.2cm}p{1.2cm}p{1.2cm}p{1.2cm}}
\toprule
{Problem} &$m$&$q$& $n$&$d_A $&$d_B$&$d_Q$  \\

\hline
  QP-D1  &800&200&2000&1&1&1\\
  QP-D2  &2000&4500&5000&1&1&1\\
  QP-D3  &100&5000&10000&1&1&1\\
  QP-D4  &4000&9000&10000&1&1&1\\
  QP-S1  &10&10&5000&0.01&0.003&9.48E-5\\
  QP-S2  &4000&10000&10000&0.01&0.001&1.10E-3\\
  QP-S3  &8000&10000&20000&0.001&0.0001&1.49E-4\\
  QP-S4  &12000&29999&30000&0.001&0.0001&3.33E-4\\
  QP-S5  &1000&4000&50000&0.001&0.0001&1.20E-4\\
  QP-S6  &1000&99000&100000&0.001&0.00002&4.39E-5\\

\toprule
\end{tabular}
\vspace{-0.4cm}
\label{table1}
\end{table}

\begin{table}[!htb]
\caption{\scriptsize{Experiments on randomly generated QPs. $f_{*}$, $f_{S}$ denote the optimal values obtained by PAL-Hom and the other corresponding solvers, ``OT'' denote the computation times (seconds)  more than 25000s.  The bolded computation times of PAL-Hom denotes that they are smaller than those of the other solvers.}}
 \tiny
 \tiny
\begin{tabular}{p{1.4cm}p{1.5cm}p{1.5cm}p{1.5cm}p{1.5cm}p{1.5cm}}
\toprule
\multirow{1}{*}{Problem} &
\multirow{3}{*}{Results}&
\multirow{3}{*}{PAL-Hom} &
\multirow{3}{*} {IPM(cplex)}&
\multirow{3}{*}{AS(matlab) } &
\multirow{3}{*}{PAS(qpOASES)}  \\
\multirow{1}{*}{m}&&&&\\
\multirow{1}{*}{n}&&&&\\

   \midrule
   \multirow{1}{*}{QP-D1 } &Time &\bf{39.30}&242.23&10314.69&2342.48\\
   \multirow{1}{*}{800 }&$\| Ax-b\|$&1.6E-11&2.9E-06&1.1E-11&1.2E-12\\
   \multirow{1}{*}{2,000} &$f_{*}-f_{S}$ &-&-2.2E-06&-1.3E-08&-2.4E-08\\
   \midrule
   \multirow{1}{*}{QP-D2 } &Time &\bf{35.45}&633.12&OT&OT\\
    \multirow{1}{*}{2,000 }  &$\| Ax-b\|$&6.9E-12&2.8E-07&-&-\\
  \multirow{1}{*}{5,000 }  &$f_{*}-f_{S}$ &-&-7.6E-04&-&-\\
   \midrule
   \multirow{1}{*}{QP-D3 } &Time &\bf{203.44}&5238.28&OT&OT\\
   \multirow{1}{*}{100 }  &$\| Ax-b\|$&1.7E-11&4.4E-07&-&-\\
   \multirow{1}{*}{10,000 }   &$f_{*}-f_{S}$ &-&-1.2E-03&-&-\\
   \midrule
   \multirow{1}{*}{QP-D4 } &Time &\bf{369.67}&6462.42&OT&OT\\
   \multirow{1}{*}{4,000 }  &$\| Ax-b\|$&1.2E-07&7.7E-07&-&-\\
    \multirow{1}{*}{10,000 }  &$f_{*}-f_{S}$ &-&-4.9E-03&-&-\\

   \midrule
   \multirow{1}{*}{QP-S1 } &Time &0.92&0.14&OT&74.22\\
   \multirow{1}{*}{10 }  &$\| Ax-b\|$&6.2E-09&8.2E-06&-&5.9E-13\\
   \multirow{1}{*}{5,000 }   &$f_{*}-f_{S}$ &-&-1.3E-09&-&1.3E-04\\

   \midrule
   \multirow{1}{*}{QP-S2 } &Time &{\bf 47.07}&126.44&OT&OT\\
    \multirow{1}{*}{4,000 } &$\| Ax-b\|$&3.5E-08&1.5E-05&-&-\\
    \multirow{1}{*}{10,000 }  &$f_{*}-f_{S}$ &-&5.7E-07&-&-\\
   \midrule
   \multirow{1}{*}{QP-S3 } &Time &208.33&128.72&OT&OT\\
   \multirow{1}{*}{8,000 }  &$\| Ax-b\|$&8.8E-08&5.8E-07&-&-\\
   \multirow{1}{*}{20,000 }    &$f_{*}-f_{S}$ &-&-1.3E-04&-&-\\
   \midrule
   \multirow{1}{*}{QP-S4 } &Time &{\bf 253.49}&5831.07&OT&OT\\
   \multirow{1}{*}{12,000 }  &$\| Ax-b\|$&3.4E-08&7.9E-07&-&-\\
   \multirow{1}{*}{30,000 }    &$f_{*}-f_{S}$ &-&-2.2E-03&-&-\\
   \midrule
   \multirow{1}{*}{QP-S5 } &Time &{\bf 520.98}&1784.71&OT&OT\\
   \multirow{1}{*}{1,000 } &$\| Ax-b\|$&1.7E-09&5.2E-08&-&-\\
    \multirow{1}{*}{50,000 }   &$f_{*}-f_{S}$ &-&2.6E-07&-&-\\
   \midrule
   \multirow{1}{*}{QP-S6 } &Time &{\bf 1757.49}&4416.68&OT&OT\\
    \multirow{1}{*}{1,000 }  &$\| Ax-b\|$&1.1E-08&4.4E-04&-&-\\
     \multirow{1}{*}{100,000 }  &$f_{*}-f_{S}$ &-&-2.1E-07&-&-\\

\toprule
\end{tabular}
\vspace{-0.4cm}
\label{table2}
\end{table}

\textbf{$\bullet$ QPs from CUTEr set.} In this part, we took convex QPs from the CUTEr test set\footnotemark[1]\footnotetext[1]{https://github.com/YimingYAN/QP-Test-Problems}, where we chose a subset of medium-scale QPs having up to 90,597 variables.
%%%%%%%%%%%%%%%%%%%%%%%%%%%%%%%%%%%%%%%%%%%%%%%%%%%%%%%%%%%%%%%%%%%%%%%%%%%%%%%%%%%%%%%%%%%%%%%%%
\begin{table}[!htb]
\caption{\scriptsize{Experiments on QPs from CUTEr test set: part I.  $f_{*}$, $f_{S}$ denote the optimal values obtained by PAL-Hom and the other corresponding solvers, ``OT'' denotes the computation times (seconds) more than 25000s. ``F" denotes that the algorithm does not converges in $10n$ iterations. The bolded computation times of PAL-Hom denotes that  they are smaller than those  of the other solvers.}}
  \tiny
\begin{tabular}{p{0.9cm}p{0.5cm}p{0.45cm}p{0.9cm}p{1.4cm}p{1cm}p{1cm}p{1cm}}
\toprule
\multirow{1}{*}{Problem} &
m&
n&
Results&
AL-Hom &
\multicolumn{1}{l} {IPM(cplex)}&
\multicolumn{1}{l}{AS(matlab)} &
\multicolumn{1}{l}{PAS(qpOASES)}  \\
   \midrule
   \multirow{3}{*}{aug2dcqp }&
   \multirow{3}{*}{10000 }&
   \multirow{3}{*}{20200 }&Time &0.65&0.41&OT&OT\\
                           &&&$\| Ax-b\|$&8.0E-13&3.9E-13&-&-\\
                           &&&$f_{*}-f_{S}$&-&-1.1E-03&-&-\\
   \midrule
   \multirow{3}{*}{aug2dqp }&
   \multirow{3}{*}{10000 }&
   \multirow{3}{*}{20200 }&Time &0.66&0.31&OT&OT\\
                           &&&$\| Ax-b\|$&9.8E-13&4.2E-13&-&-\\
                           &&&$f_{*}-f_{S}$&-&-1.0E-04&-&-\\
   \midrule
   \multirow{3}{*}{aug3dcqp }&
   \multirow{3}{*}{1000 }&
   \multirow{3}{*}{3873}&Time &0.14&0.04&569.22&2763.23\\
                           &&&$\| Ax-b\|$&2.4E-14&2.4E-13&2.4E-13&3.0E-15\\
                           &&&$f_{*}-f_{S}$&-&-2.7E-06&-1.4E-09&4.3E-12\\
   \midrule
   \multirow{3}{*}{aug3dqp }&
   \multirow{3}{*}{1000 }&
   \multirow{3}{*}{3873 }&Time &0.20&0.05&OT&3513.54\\
                           &&&$\| Ax-b\|$&1.9E-13&1.1E-14&-&7.4E-11\\
                           &&&$f_{*}-f_{S}$&-&-7.5E-08&-&3.3E-10\\
   \midrule
   \multirow{3}{*}{cont-050 }&
   \multirow{3}{*}{2401}&
   \multirow{3}{*}{2597  }&Time &0.92&0.27&17.92&5397.38\\
                           &&&$\| Ax-b\|$&1.5E-13&3.8E-14&1.2E-13&4.7E-14\\
                           &&&$f_{*}-f_{S}$&-&-7.2E-10&-3.3E-13&2.1E-09\\
   \midrule
   \multirow{3}{*}{cont-100 }&
   \multirow{3}{*}{9801 }&
   \multirow{3}{*}{10197 }&Time &2.18&0.56&817.25&OT\\
                           &&&$\| Ax-b\|$&4.0E-13&7.4E-14&1.7E-12&-\\
                           &&&$f_{*}-f_{S}$&-&4.1E-08&-3.3E-13&-\\
   \midrule
   \multirow{3}{*}{cont-101 }&
   \multirow{3}{*}{10098}&
   \multirow{3}{*}{10197}&Time &4.18&0.84&848.61&OT\\
                           &&&$\| Ax-b\|$&7.3E-13&8.5E-10&3.3E-12&-\\
                           &&&$f_{*}-f_{S}$&-&4.2E-07&-2.4E-06&-\\
      \midrule
   \multirow{3}{*}{cont-200 }&
   \multirow{3}{*}{39601}&
   \multirow{3}{*}{40397}&Time &11.00&1.40&OT&OT\\
                           &&&$\| Ax-b\|$&5.5E-13&1.5E-13&-&-\\
                           &&&$f_{*}-f_{S}$&-&7.3E-07&-&-\\
      \midrule
   \multirow{3}{*}{cont-201 }&
   \multirow{3}{*}{40198}&
   \multirow{3}{*}{40397}&Time &23.13&2.31&OT&OT\\
                           &&&$\| Ax-b\|$&7.2E-10&1.8E-08&-&-\\
                           &&&$f_{*}-f_{S}$&-&1.7E-06&-&-\\
      \midrule
   \multirow{3}{*}{cont-300 }&
   \multirow{3}{*}{90298}&
   \multirow{3}{*}{90597}&Time &41.46&4.62&OT&OT\\
                           &&&$\| Ax-b\|$&8.4E-09&2.8E-08&-&-\\
                           &&&$f_{*}-f_{S}$&-&4.1E-05&-&-\\
   \midrule
   \multirow{3}{*}{cvxqp1$\_$l }&
   \multirow{3}{*}{5000}&
   \multirow{3}{*}{10000}&Time &50.22&24.65&OT&OT\\
                           &&&$\| Ax-b\|$&9.9E-08&5.4E-07&-&-\\
                             &&&$f_{*}-f_{S}$&-&-2.5E-02&-&-\\
   \midrule
   \multirow{3}{*}{cvxqp1$\_$m }&
   \multirow{3}{*}{500 }&
   \multirow{3}{*}{1000 }&Time &\bf{0.48}&0.87&7.61 &92.43\\
                           &&&$\| Ax-b\|$&7.2E-13&1.9E-07&8.2E-14 &2.5E-13\\
                           &&&$f_{*}-f_{S}$&-&-1.7E-04&1.2E-05 &-1.6E-04\\
   \midrule
   \multirow{3}{*}{cvxqp1$\_$s }&
   \multirow{3}{*}{50 }&
   \multirow{3}{*}{100 }&Time &{\bf 0.01}&0.02& 0.02&0.15\\
                           &&&$\| Ax-b\|$&2.1E-14&9.5E-12&9.8E-15 &2.9E-15\\
                           &&&$f_{*}-f_{S}$&-&-1.5E-05&3.1E-09 &-5.2E-07\\
   \midrule
   \multirow{3}{*}{cvxqp2$\_$l }&
   \multirow{3}{*}{2500 }&
   \multirow{3}{*}{10000 }&Time &\bf{1.78}&12.63&OT&OT\\
                           &&&$\| Ax-b\|$&1.8E-08&1.2E-08&-&-\\
                           &&&$f_{*}-f_{S}$&-&-1.9E-03&-&-\\
   \midrule
   \multirow{3}{*}{cvxqp2$\_$m }&
   \multirow{3}{*}{250 }&
   \multirow{3}{*}{1000}&Time &\bf{0.15}&0.53&21.72&26.37\\
                           &&&$\| Ax-b\|$&7.5E-08&3.8E-08& 3.8E-14&7.7E-15\\
                           &&&$f_{*}-f_{S}$&-&-1.3E-03&4.3E-06&-3.1E-07\\
   \midrule
   \multirow{3}{*}{cvxqp2$\_$s }&
   \multirow{3}{*}{25 }&
   \multirow{3}{*}{100}&Time &\bf{0.01}&0.02&0.04 &0.05\\
                           &&&$\| Ax-b\|$&6.8E-08&3.8E-10&8.5E-15 &2.6E-15\\
                           &&&$f_{*}-f_{S}$&-&-4.0E-05& -5.0E-08&1.3E-06\\
   \midrule
   \multirow{3}{*}{cvxqp3$\_$l }&
   \multirow{3}{*}{7500 }&
   \multirow{3}{*}{10000 }&Time &54.44&30.19&OT&OT\\
                           &&&$\| Ax-b\|$&4.8E-05&1.9E-05&- &- \\
                           &&&$f_{*}-f_{S}$&-&-6.7E-04&-&-\\
   \midrule
   \multirow{3}{*}{cvxqp3$\_$m }&
   \multirow{3}{*}{750}&
   \multirow{3}{*}{1000  }&Time &\bf{0.66}&0.95&10.11 &292.15\\
                           &&&$\| Ax-b\|$&7.5E-08&3.6E-09& 1.7E-13&4.4E-12\\
                           &&&$f_{*}-f_{S}$&-&-5.6E-02& 2.3E-05&1.2E-04\\
   \midrule
   \multirow{3}{*}{cvxqp3$\_$s }&
   \multirow{3}{*}{75}&
   \multirow{3}{*}{100  }&Time &{\bf 0.01}&0.01&0.03 &0.22\\
                           &&&$\| Ax-b\|$&1.8E-09&6.6E-12&1.7E-14 &1.9E-14\\
                           &&&$f_{*}-f_{S}$&-&-9.7E-06& 2.7E-08&-2.1E-07\\

 \midrule
   \multirow{3}{*}{gouldqp2 }&
   \multirow{3}{*}{349}&
   \multirow{3}{*}{699 }&Time &\bf{0.00}&0.01&0.11&0.02\\
                           &&&$\| Ax-b\|$&0.0E+00&1.9E-08&0.0E+00&0.0E+00\\
                           &&&$f_{*}-f_{S}$&-&-3.0E-11&0.0E+00&0.0E+00\\

   \midrule
   \multirow{3}{*}{gouldqp3 }&
   \multirow{3}{*}{349 }&
   \multirow{3}{*}{699}&Time &0.04&0.02&0.15&55.74\\
                           &&&$\| Ax-b\|$&1.3E-11&3.8E-09&1.0E-13&1.5E-14\\
                           &&&$f_{*}-f_{S}$&-&-3.1E-06&2.3E-06&-2.1E-10\\

\toprule
\end{tabular}
\vspace{-0.6cm}
\label{table3}
\end{table}
%%%%%%%%%%%%%%%%%%%%%%%%%%%%%%%%%%%%%%%%%%%%%%%%%%%%%%%%%%%%%%%%%%%%%%%
\begin{table}[!htb]
\caption{\scriptsize{Experiments on QPs from CUTEr test set: part II. $f_{*}$, $f_{S}$ denote the optimal values obtained by PAL-Hom and the other corresponding solvers, ``OT'' denotes the computation times (seconds) more than 25000s.  ``F" denotes that the algorithm does not converges in $10n$ iterations. The bolded computation times of PAL-Hom denote that  they are smaller than those  of the other solvers.}}
 \tiny
 \tiny
\begin{tabular}{p{0.9cm}p{0.5cm}p{0.45cm}p{0.9cm}p{1.4cm}p{1cm}p{1cm}p{1cm}}
\toprule
\multirow{1}{*}{Problem} &
m&
n&
Results&
PAL-Hom &
\multicolumn{1}{l} {IPM(cplex)}&
\multicolumn{1}{l}{AS(matlab)} &
\multicolumn{1}{l}{PAS(qpOASES)}  \\

   \midrule
   \multirow{3}{*}{powell20 }&
   \multirow{3}{*}{10000 }&
   \multirow{3}{*}{10000 }&Time &\bf{0.14}&0.17&OT&OT\\
                           &&&$\| Ax-b\|$&2.5E-08&1.3E-05&-&-\\
                           &&&$f_{*}-f_{S}$&-&3.9E-03&-&-\\
   \midrule
   \multirow{3}{*}{qgrow7 }&
  \multirow{3}{*}{140 }&
   \multirow{3}{*}{301 }&Time &0.37&0.01&F&1.56\\
                           &&&$\| Ax-b\|$&2.5E-08&1.0E-10&-&1.6E-02\\
                           &&&$f_{*}-f_{S}$&-&-2.3E-03&-&1.4E-01\\
   \midrule
   \multirow{3}{*}{qgrow15 }&
    \multirow{3}{*}{300 }&
   \multirow{3}{*}{645 }&Time &1.55&0.02&F&1936.40\\
                           &&&$\| Ax-b\|$&2.8E-08&2.5E-10&-&1.5E-14\\
                           &&&$f_{*}-f_{S}$&-&-8.6E-03&-&-4.2E-11\\
   \midrule
   \multirow{3}{*}{qgrow22 }&
   \multirow{3}{*}{440 }&
   \multirow{3}{*}{946 }&Time & 0.57&0.28&F&5924.02\\
                           &&&$\| Ax-b\|$&6.1E-11&4.6E-07&-&4.9E-12\\
                           &&&$f_{*}-f_{S}$&-&7.2E-05&-&-3.3E-09\\
   \midrule
   \multirow{3}{*}{qscsd1 }&
   \multirow{3}{*}{77 }&
   \multirow{3}{*}{760 }&Time &0.10&0.02&F&1.44\\
                           &&&$\| Ax-b\|$&4.0E-08&3.7E-11&-&9.1E-10\\
                           &&&$f_{*}-f_{S}$&-&-1.1E-08&-&3.2E-08\\
   \midrule
   \multirow{3}{*}{qscsd6}&
   \multirow{3}{*}{147 }&
   \multirow{3}{*}{1350 }&Time &0.17&0.02&F&11.26\\
                           &&&$\| Ax-b\|$&3.9E-08&2.1E-12&-&9.1E-10\\
                           &&&$f_{*}-f_{S}$&-&-5.1E-08&-&3.3E-09\\
   \midrule
   \multirow{3}{*}{qscsd8 }&
   \multirow{3}{*}{397 }&
   \multirow{3}{*}{2750 }&Time &0.89&0.03&F&36.87\\
                           &&&$\| Ax-b\|$&3.0E-08&2.5E-08&-&2.5E-12\\
                           &&&$f_{*}-f_{S}$&-&-1.0E-06&-&3.2E-09\\
   \midrule
   \multirow{3}{*}{stcqp1 }&
   \multirow{3}{*}{2052}&
   \multirow{3}{*}{4097 }&Time &0.13&0.07&1357.24&499.14\\
                           &&&$\| Ax-b\|$&4.0E-09&4.0E-13&4.1E-12&8.1E-13\\
                           &&&$f_{*}-f_{S}$&-&-1.3E-04&3.3E-06&-3.7E-08\\
   \midrule
   \multirow{3}{*}{stcqp2}&
   \multirow{3}{*}{2052 }&
   \multirow{3}{*}{4097 }&Time &{\bf0.06}&0.67&916.18&1420.90\\
                           &&&$\| Ax-b\|$&3.1E-08&0.0E+00&3.6E-11&3.7E-11\\
                           &&&$f_{*}-f_{S}$&-&-1.1E-04&-3.7E-11&2.2E-08\\

\toprule
\end{tabular}

\label{table4}
\end{table}
%%%%%%%%%%%%%%%%%%%%%%%%%%%%%%%%%%%%%%%%%%%%%%%%%%%%%%%

We compared PAL-Hom with the IPM solver in CPLEX 12.6, the PAS solver in qpOASES and the AS solver in MATLAB 2013a. The comparison includes the computation time (seconds), equality constraint violations and  optimal values. The results are reported in Tables \ref{table2}-\ref{table4}.

The numerical results  show that PAL-Hom is effective at solving these QPs. PAL-Hom outperforms the AS solver in MATLAB and the PAS solver in qpOASES  and is competitive with the IPM solver in CPLEX, especially for the randomly generated problems.

Moreover,  to show that the homotopy algorithm with warm start by APG is  meaningful for the augmented Lagrangian subproblems, we compared   the algorithm with the IPM solver in CPLEX and Hager et al.'s  active-set algorithm (ASA) \cite{hager2006new}, which consists of a nonmonotone gradient projection step, an unconstrained optimization step, and a set of rules for branching between the two steps. ASA is shown to be faster than TRON \cite{lin1999newton} for solving the 50 box-constrained problems in the CUTEr library \cite{bongartz1995cute} and competitive with TRON for the 23 box-constrained problems in the
MINPACK-2 library \cite{averick1992minpack}. Furthermore,  to show that the homotopy tracking with the sorting technique and the $\varepsilon$-precision verification and correction technique  is  more efficient than the PAS solver in qpOASES for solving parametric nonnegative QP problems, we compared  it with the  PAS solver in qpOASES for solving the first augmented Lagrangian subproblem (\ref{almsub}) from $\hat{z}$.

The results are reported in Table \ref{table5}. Clearly, the homotopy tracking is much faster than the PAS solver in qpOASES. Moreover, from the results, we see that APG is efficient at predicting the optimal active set; that is, from the approximate solution, a small number of  tracking steps is required to obtain an exact solution. In addition, we see that the homotopy algorithm is robust for   problems with large condition numbers. Simultaneously, the results demonstrate that the homotopy algorithm is clearly faster than PAS(with initial point $\hat{z}$), ASA and IPM for solving the augmented Lagrangian subproblems.

\begin{table}[!htb]
\caption{\scriptsize{IPM(CPLEX), ASA, PAS solver in qpOASES and the homotopy algorithm solving the first augmented Lagrangian problem. ``Total" denotes the computation times (seconds) of APG and the homotopy tracking steps together, ``Hom-tra." denotes the computation times of the homotopy tracking steps starting from $\hat{z}$. The computation times of PAS method denote the cost of solving the parametric quadratic programming from $\hat{z}$. ``OT'' denotes the computation times more than 25000s. $C_H$ denotes the condition number of $H$. The bolded computation times (APG+Hom-tra.) in the ``Total" column  denotes that they are smaller than those of the other solvers, and the bolded computation times  in the ``Hom-tra." column  denotes that  they are smaller than that  of the PAS solver.}}

\vspace{-0.1cm}
\begin{center} \tiny
\begin{tabular}{llllllllllllll}
\toprule
\multirow {2}{*}{Problem}&
\multirow {2}{*}{n}&
\multirow {2}{*}{$C_H$}&
\multicolumn{3}{c}{Homotopy}&
&
\multicolumn{2}{c}{PAS(qpOASES)}&ASA&IPM(cplex)

  \\
\cline{4-6}
\cline{8-9}
   &&&Total&Hom-tra. &Iter.& & Time &Iter.&Time &Time \\
\midrule
aug2dcqp   &20200&5.2E+06&\bf{0.44}&\bf{0.16}&9&&221.33&9&1.90&4.24 \\
aug2dqp   &20200&8.0E+13&\bf{0.43}&\bf{0.15}&8&&273.11&10&12.19&4.30\\
aug3dcqp   &3873&8.6E+03&\bf{0.06}&\bf{0.02}&2&&16.33&2&0.07& 2.01\\
aug3dqp   &3873&2.8E+11&\bf{0.07}&\bf{0.02}&2&&11.44&2&0.12& 1.96\\
cont-50   &2597&3.4E+08&\bf{0.48}&\bf{0.11}&4&&1.61&4 &45.57&0.56\\
cont-100   &10197&1.1E+08&\bf{0.84}&\bf{0.21}&4&&332.13&4&15.50&4.41\\
cont-101  &10197&2.6E+10&\bf{2.22}&\bf{0.62}&21&&61.21&25&20990.52&4.55 \\
cont-200  &40397&1.2E+10&\bf{4.46}&\bf{0.35}&5&&171.33&8&OT&27.87 \\
cont-201  &40397&4.5E+13&\bf{12.88}&\bf{1.77}&38&&864.77&44&OT&35.45 \\
cont-300  &90597&8.0E+10&\bf{27.99}&\bf{4.99}&66&&OT&OT&OT&316.99 \\
cvxqp1$\_$l  &10000&8.8E+13&\bf{5.99}&\bf{1.66}&45&&277.33&49&406.54&42.11 \\
cvxqp2$\_$l  &10000&6.2E+10&\bf{0.56}&\bf{0.13}&10&&162.11 &11&7.10&40.74\\
cvxqp3$\_$l  &10000&6.3E+08&\bf{8.33}&\bf{2.16}&25&& 311.23&29&114.80&46.22 \\
powell20  &10000&4.0E+08&\bf{0.04}&\bf{0.01}&2&&0.10&4&163.58&0.12 \\
qgrow22  &946&6.3E+06&\bf{0.30}&\bf{0.03}&3&&0.16&3&0.03&0.41 \\
qscsd6  &1350&3.0E+10&\bf{0.11}&\bf{0.03}&33&&0.33&37&14459.72& 0.75\\
qscsd8  &2740&1.9E+10&\bf{0.40}&\bf{0.15}&72&&3.27&81&351.74&0.65 \\
stcqp1  &4097&1.3E+03&0.04&\bf{0.01}&2&&0.38&2&0.04&1.27\\
stcqp2  &4097&8.2E+02&\bf{0.03}&\bf{0.01}&1&&5.43&1&0.04&1.72 \\
   \bottomrule
\end{tabular}
\end{center}
\vspace{-0.6cm}
\label{table5}
\end{table}

Moreover,  to show the adaptability of  PAL-Hom for the degenerate QP problems, we tested the homotopy method on solving highly ill-conditioned non-negative constrained QP problems. We  conducted experiments like this for the efficiency of PAL-Hom depends on the solving of  the proximal augmented Lagrangian subproblems which are degenerate non-negative constrained QP problems with proximal terms. We generated the non-negative constrained QP problems (\ref{almsub}) with MATLAB codes  as follows.

$d={\rm zeros}(n,1);~d(1:\gamma n)=\frac{L_{\max}}{\gamma n}(1:\gamma n);~H=U'*{\rm diag}(d)*U+\delta*I$;

$f={\rm randn}(n,1)$;

\noindent where $0<\gamma<1$ denotes the ratio of nonzero eigen-values, $U$ is the unitary matrix.

We generated the problems with $n=2000, L_{\max}=10^8$ and $\gamma\in\{0.1,0.5\}$. We adjusted the condition number of $H$ which equals to $\frac{L_{\max}+\delta}{\delta}$ by changing the value of $\delta$. $H$ has $\gamma n$ egien-values bigger than $\delta$ and the rest eigen-values are $\delta$.
Obviously, $H$ is ill-conditioned when $\delta$ is small.

\begin{figure}[!htb]
  \centering

  \includegraphics[width=5.5cm]{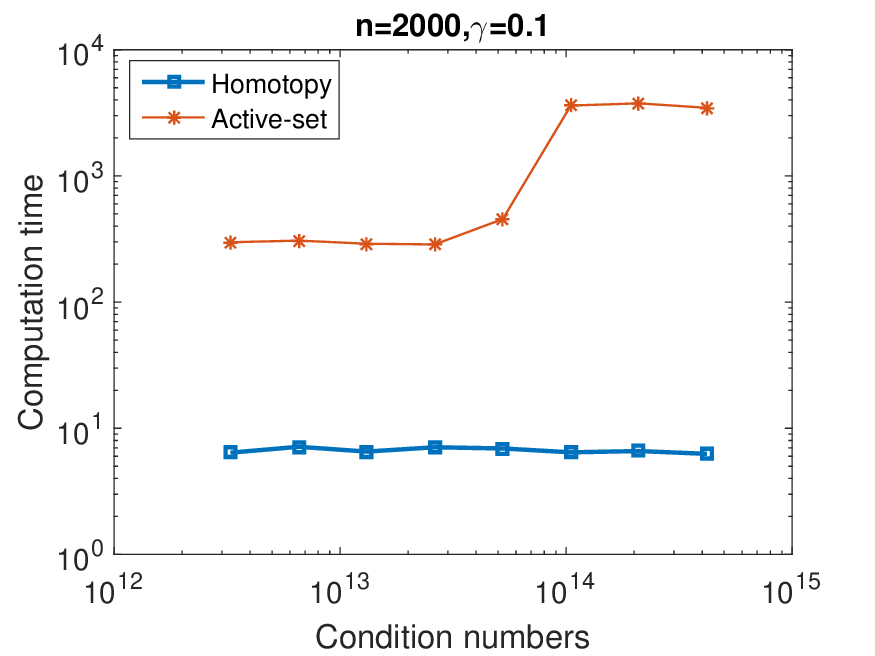}
  \includegraphics[width=5.5cm]{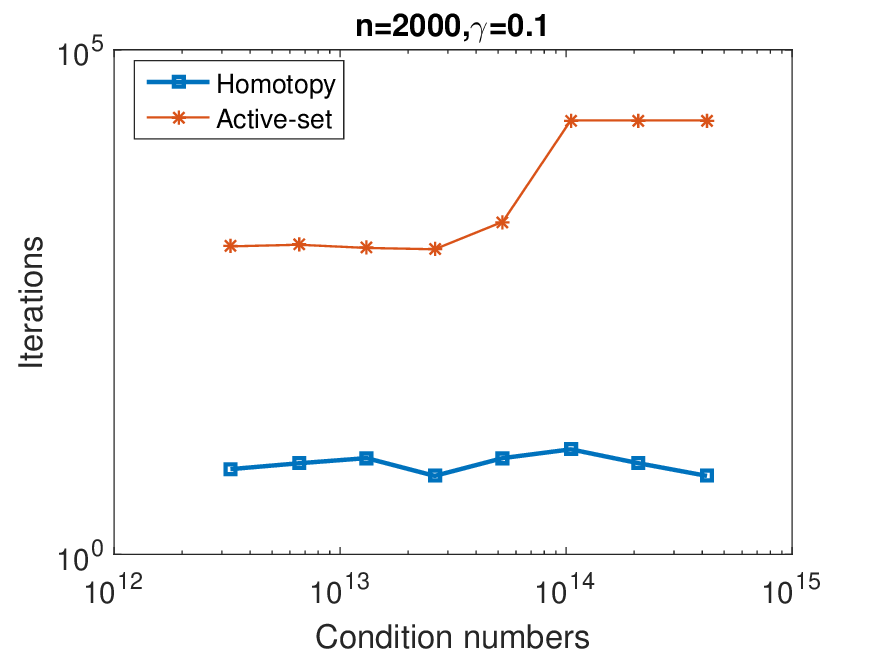}\\
    \includegraphics[width=5.5cm]{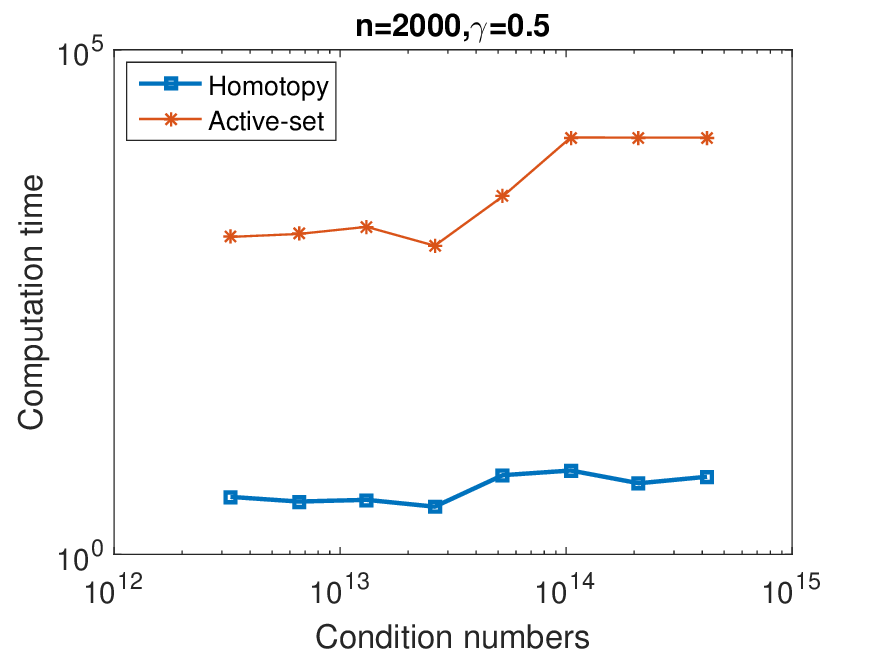}
  \includegraphics[width=5.5cm]{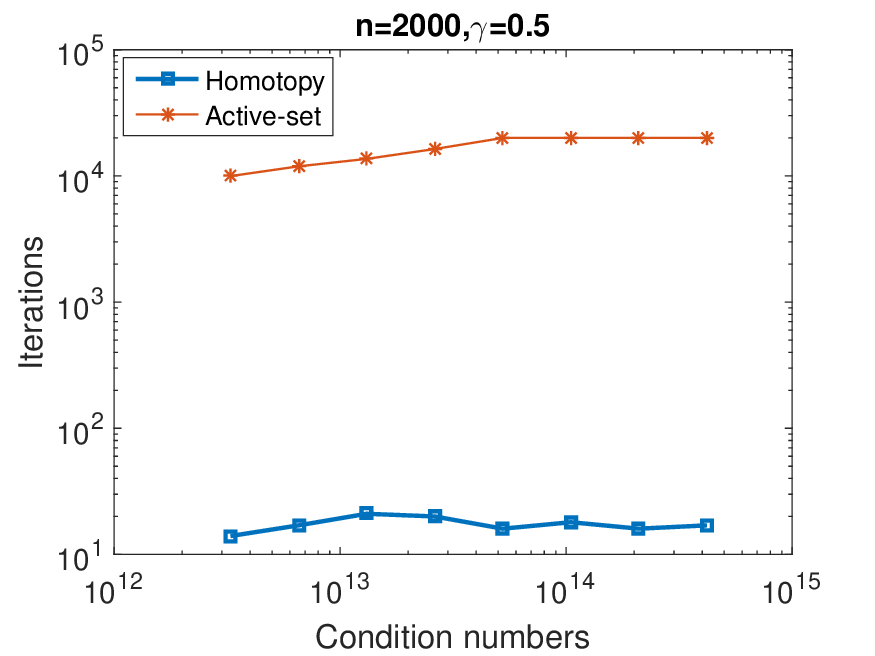}\\
  \caption{The homotopy  and the active-set method for ill-conditioned non-negative constrained QP problems. The results include the computation time (seconds) and  the number of the iterations.}\label{figure2}
\end{figure}

We solved these problems by the homotopy algorithm and the active-set method, respectively. The maximum iterations of the active-set method is set to $100*\gamma*n$. The results are shown in Figure \ref{figure2}.

The results demonstrate that the homotopy method is robust for the ill-conditioned non-negative constrained QP problems, while the active-set method requires much more time  when the condition number is larger. Moreover, the number of steps of the homotopy tracking does not change much when the condition number increases,  while the AS method often exceed the maximum iterations when the condition number is very large.

\subsection{QPs from SVM for recognition.} In this section, we tested AL-Hom for solving QPs from SVMs that were applied to handwritten digit recognition and speech recognition. Given the training set $\{X_i,y_i\}_{i=1}^{n}$ and testing set $\{T_j,s_j\}_{j=1}^{n_1}$, where $X_i, T_j$ are feature vectors and $y_i, s_j\in\{-1,+1\}$ are the labels, the SVM classifies the testing set by a classifier
$$f(x)={\rm sign}(\sum_{i=1}^{n}y_i\alpha_i^{*}K(X_i,x)+b^*),$$
where $K$ is called the kernel function, $b^*=y_j-\sum_{i=1}^{n}y_i\alpha_i^{*}K(X_i,X_j)$, for some $ \alpha_{j}^{*}>0$, and $\alpha^*$ is the solution of the following problem
\begin{eqnarray} \nonumber
&&\min_{\alpha} ~~~~\frac{1}{2}\sum_{i=1}^{n}\sum_{i=1}^{n}y_iy_j\alpha_i\alpha_j K(X_i,X_j)-\sum_{i=1}^{n}\alpha_i\\\label{svm}
&&{\rm s.t.}~~~~\sum_{i=1}^{n}y_i\alpha_i=0,\\\nonumber
&& ~~~~~~~~~0\leq\alpha_i\leq C, i=1,...,n,
\end{eqnarray}
which is a dense QP problem.

\begin{figure}[!htb]
  \centering
  \vspace{-0cm}
  \includegraphics[width=9.5cm]{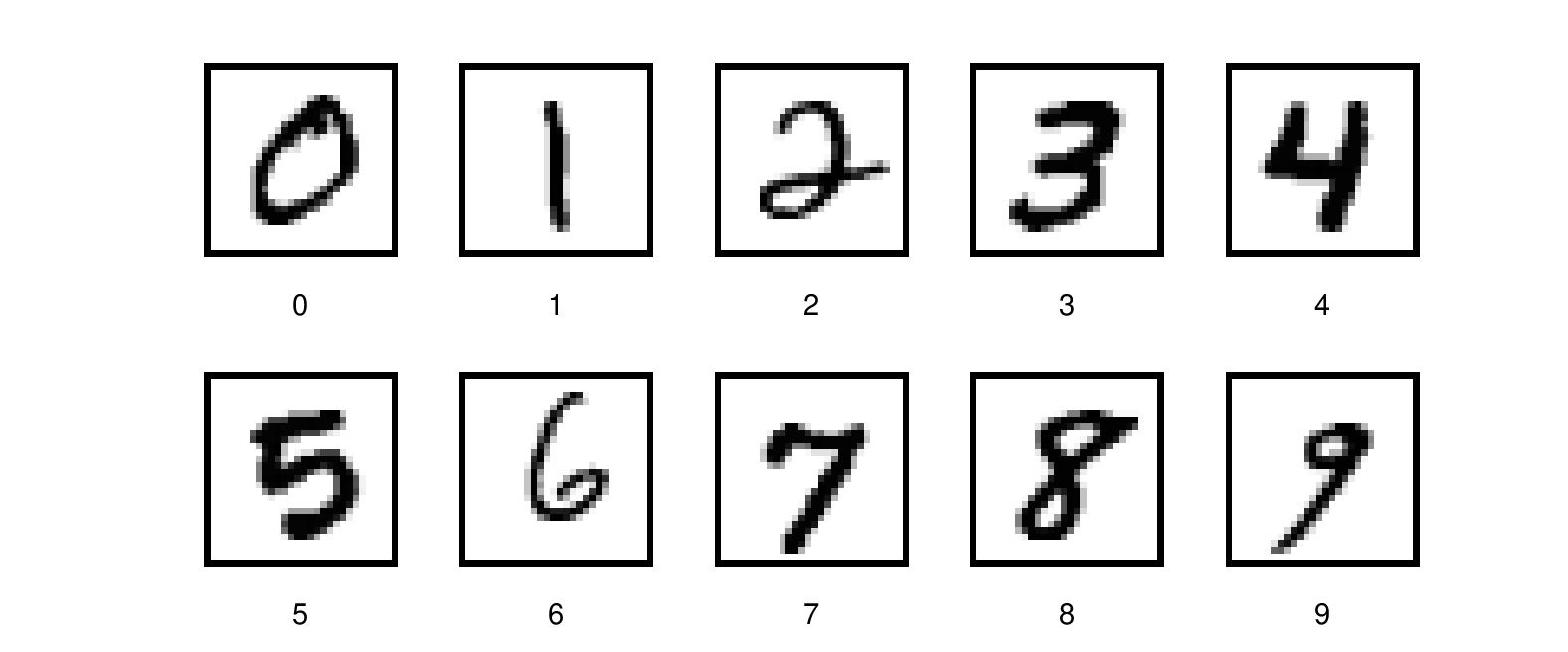}\\
  \vspace{-0cm}
  \caption{Ten classifications of the MNIST handwritten digit database}
\label{figure3}
\end{figure}

We conducted the experiments with three databases. The first database is the isolated letter speech database from UCI\cite{Lichman2013}, which contains a training set with 6,238 samples and a testing set with 1,559 samples. This database has 26 classifications: i.e., A-Z, and every sample has 617 attributes. The second one is   the MNIST database of handwritten digits\footnotemark[2]\footnotetext[2]{http://yann.lecun.com/exdb/mnist}, which contains a training set with 60,000 samples and a testing set with 10,000 samples. Every sample is one $28\times 28$ pixel picture, that is, every sample has 784 attributes. This database has ten classifications as shown in Figure \ref{figure3}. The third database is the web page classification task, which is included in the LIBSVM database set\footnotemark[3]\footnotetext[3]{https://www.csie.ntu.edu.tw/~cjlin/libsvmtools/datasets/}. This database contains 8 training sets and 8 testing sets of different sizes. Every set contains samples divided into 2 classes and every sample has 300 features. Because the training sets have repetitive samples, we processed them individually by removing the repetitive samples.

($\ref{svm}$) is a  model for the 2-class classification problem; however, the letter speech and MNIST databases are multiclassification problems. Therefore. we handled the multiclassification problems with two strategies.

 The first strategy is that  for any classification $p, p=1,..,P$, where $P$ denotes the number of   classifications, we obtained $\alpha^{*,p}, b^{*,p}$ by solving ($\ref{svm}$) with
\begin{eqnarray} \nonumber
y^p=\left\{
\begin{array}{lcl}
1,      &      & X_i\in cl.p;\\
-1,  &      & else,
\end{array} \right.
\end{eqnarray}
where $cl.p$ denotes the $p$-th classification. Then, we have the first classifier for multiclassification problems as follows
\begin{eqnarray} \label{classifier}
f_1(x)=\arg\max_{p}~ (y^p_i\alpha_i^{*,p}K(X_i,x)+b^{*,p}).
\end{eqnarray}

The second strategy is that for any $p\neq q\in\{1,...,P\}$, choose the samples from the training set whose corresponding labels are $p$ or $q$, and let
$$y^{p,q}=\left\{
\begin{array}{lcl}
1,      &      & X_i\in cl.p;\\
-1,  &      &  X_i\in cl.q,
\end{array} \right.$$
then, we  have the second  classifier for multiclassification problems  as follows.
\begin{eqnarray} \label{classifier2}
f_2(x)=\arg\max_{p}~ \sum_{q} {\rm sign} (y^{p,q}_i\alpha_i^{*,p,q}K(X_i,x)+b^{*,p,q}).
\end{eqnarray}

The first strategy needs to solve $(\ref{svm})$  $P$ times, and the size of every problem is equal to the number of   samples in the training set. The second strategy needs to solve $(\ref{svm})$ $\frac{P(P-1)}{2}$ times, however, it only needs to solve a problem of size equal to the number of  samples of the classification $p$ and $q$ together each time. In LIBSVM, the  multiclassification classifier adopts the second strategy.

In our experiments, we used the polynomial kernel
\begin{eqnarray} \nonumber
K(x,y)=(x\cdot y+c)^d
\end{eqnarray}
with $c=0$ and  $d=2$ for the spoken letter database and the MNIST databases and  the Gaussian kernel
\begin{eqnarray} \nonumber
K(x,y)=e^{-\sigma\|x-y\|^2}
\end{eqnarray}
with $\sigma=0.1$ for the web pages classification task.

Because the QP problems in this section are strictly convex, we compared AL-Hom with the IPM solver in CPLEX and the sequential minimal optimization (SMO) method \cite{fan2005working} in LIBSVM 3.22 \cite{CC01a} which is a well-known package for SVMs. We report the results in Tables \ref{table6} and \ref{table7}, where ``Err.1'' and ``Err.2'', respectively, denote the number of  misclassifications of classifier $(\ref{classifier})$ and classifier $(\ref{classifier2})$ for the test set. Here, we simply list the computation time of the first strategy and not that of the second strategy because it contains $\frac{P(P-1)}{2}$ parts.  We only give the total time of the second strategy in the title of the tables. Moreover, we do not list the results of IPM for the MNIST database because it took substantially more time than the other two algorithms. The results show that the IPM solver in CPLEX faces difficulties in solving the QPs from SVMs. PAL-Hom outperforms IPM. Moreover, although AL-Hom does not exploit the structure or particularity, it is competitive with SMO which extensively exploits the structure and particularity of the SVM problem. We believe that AL-Hom would be more competitive for SVMs if we implement AL-Hom   using the structure and the particularity the SVM problem, such as utilizing the framework of Osuna's decomposition algorithm  \cite{osuna1997improved}.

\begin{table}[!htb]
 \caption{\scriptsize{Experiments on QPs from SVM  for the classification of the isolated letter speech database. The computation time (seconds) in this table is the cost  of the first strategy. In the second strategy, AL-Hom took 4.92s, IPM took 107.08s and SMO took 13.37s in all. The bolded computation times of AL-Hom denotes that they are  smaller than those of the other solvers.}}

    \label{table6}
    \tiny

    \begin{tabular}{llllp{0.3cm}lllp{0.3cm}llllllllllll}

    \toprule
    \multirow{2}{*}{Class.}& \multicolumn{3}{c}{AL-Hom} && \multicolumn{3}{c}{IPM (cplex)}&& \multicolumn{3}{c}{SMO (LIBSVM)} \\
    \cline{2-4}\cline{6-8}\cline{10-12}
    &Time/s&Err.1&Err.2&&Time/s&Err.1&Err.2&&Time/s&Err.1&Err.2\\
    \hline
    cl.A &\bf{3.43}&0&0&&240.80&0&0&&6.65&0&0\\
    cl.B &\bf{4.33}&5&4&&242.58&5&4&&7.23&5&4\\
    cl.C &\bf{1.77}&0&0&&259.33&0&0&&5.22&0&0\\
    cl.D &\bf{5.01}&4&3&&233.34&4&3&&6.86&4&3\\
    cl.E &\bf{3.10}&0&2&&255.32&0&2&&6.77&0&2\\
    cl.F &\bf{3.16}&0&2&&232.10&0&2&&6.84&0&2\\
    cl.G &\bf{3.95}&0&0&&260.34&0&0&&6.98&0&0\\
    cl.H &\bf{1.88}&0&0&&220.50&0&0&&4.82&0&0\\
    cl.I &\bf{1.84}&1&1&&223.03&1&1&&5.37&1&1\\
    cl.J &\bf{1.80}&1&1&&273.58&1&1&&6.45&1&1\\
    cl.K &\bf{3.87}&2&2&&218.09&2&2&&6.94&1&2\\
    cl.L &\bf{1.74}&0&0&&200.01&0&0&&5.33&0&0\\
    cl.M &\bf{3.66}&9&7&&252.48&9&6&&5.34&9&6\\
    cl.N &\bf{4.32}&8&9&&230.56&8&9&&6.73&8&9\\
    cl.O &\bf{3.20}&0&0&&208.38&0&0&&7.28&0&0\\
    cl.P &6.73&0&6&&235.25&0&5&&5.41&0&5\\
    cl.Q &\bf{1.71}&4&0&&226.73&4&0&&7.78&4&0\\
    cl.R &\bf{1.44}&0&0&&258.80&0&0&&5.22&0&0\\
    cl.S &\bf{1.54}&3&3&&229.37&3&3&&4.94&3&3\\
    cl.T &\bf{4.63}&3&5&&231.48&3&6&&7.30&3&6\\
    cl.U &\bf{1.90}&2&2&&225.77&2&2&&5.88&2&2\\
    cl.V &\bf{5.33}&5&6&&231.16&5&5&&7.02&5&5\\
    cl.W &\bf{2.38}&0&0&&256.07&0&1&&6.57&0&1\\
    cl.X &\bf{1.63}&0&0&&228.58&0&0&&5.02&0&0\\
    cl.Y &\bf{1.44}&0&0&&220.33&0&0&&4.49&0&0\\
    cl.Z &\bf{2.13}&4&3&&247.07&4&3&&5.91&4&3\\
    Total &\bf{77.93}&51&56&&6138.86&51&55&&161.59&51&55\\

    \bottomrule
    \end{tabular}

\end{table}

\begin{table}[!htb]
 \caption{\scriptsize{Experiments on QPs from SVM for the classification of the mnist database.  The computation time (seconds) in this table is the cost  of the first strategy. In the second strategy, AL-Hom took 394.33s and SMO took 282.43s in all. The bolded computation times of AL-Hom denotes that they are smaller than those of the other solvers.}}

    \label{table7}
    \ \tiny

    \begin{tabular}{p{1cm}p{1cm}p{1cm}p{1cm}p{1.5cm}p{1cm}p{1cm}p{1cm}p{1cm}p{1cm}p{1cm}p{1cm}p{1cm}p{1cm}}
    \toprule
    \multirow{2}{*}{Class.}& \multicolumn{3}{c}{AL-Hom} && \multicolumn{3}{c}{SMO (LIBSVM)} \\
    \cline{2-4}\cline{6-8}
    &Time/s&Err.1&Err.2&&Time/s&Err.1&Err.2\\
    \hline
    cl.0 &\bf{901.33}&10&7 &&939.29&10&8\\
    cl.1 &644.27&10&7 &&540.36&10&8\\
    cl.2 &\bf{2001.13}&22&24 &&2336.39&22&24\\
    cl.3 &\bf{2225.65}&23& 23&&3501.45&23&25\\
    cl.4 &2001.42&17& 16&&1492.39&17&16\\
    cl.5 &\bf{1978.43}&19&20 &&2397.18&19&19\\
    cl.6 &1117.33&17&18 &&1057.11&17&18\\
    cl.7 &\bf{1863.33}&22&25 &&1926.60&22&26\\
    cl.8 &\bf{2854.22}&24&22 &&4028.26&24&23\\
    cl.9 &\bf{2131.33}&29& 30&&4111.90&29&28\\
    Total. &\bf{17718.44}&193&192&&22230.8&193&195\\
    \bottomrule
    \end{tabular}

\end{table}

\begin{table}[!htb]
 \caption{\scriptsize{Experiments on QPs from SVM for the classification of the web classification task. ``OT'' denotes the computation times more than 25000s.}}
    \label{table8}
    \tiny

    \begin{tabular}{lllllllllllllllllllllll}
    \toprule
    \multirow{2}{*}{Problem}& \multirow{2}{*}{Training.set}& \multirow{2}{*}{Testing.set}& \multicolumn{2}{c}{AL-Hom} && \multicolumn{2}{c}{IPM (cplex)}&& \multicolumn{2}{c}{SMO (LIBSVM)} \\
    \cline{4-5}\cline{7-8}\cline{10-11}
    &&&Time/s&Err.&&Time/s&Err.&&Time/s&Err.\\
    \hline
    w1a &2123&47272&0.93&1030&&6.40&1030&&0.40&1030\\
    w2a &2950&46279&1.37&895&&18.01&895&&0.64&895\\
    w3a &4108&44837&2.65&823&&54.81&823&&1.02&823\\
    w4a &6049&42383&5.34&728&&205.78&728&&1.88&728\\
    w5a &7970&39861&11.22&648&&620.39&648&&3.02&648\\
    w6a &13268&32561&34.78&419&&2950.54&419&&7.14&419\\
    w7a &18530&25057&103.13&321&&8766.35&321&&25.64&321\\
    w8a &34704&14951&553.74&113&&OT&-&&165.54&113\\
    \bottomrule
    \end{tabular}

\end{table}

AL-Hom is much faster than IPM implemented in CPLEX for solving QPs from SVMs for the following two reasons: first, ALM is effective for SVM optimization because it requires only several iterations to achieve a satisfactory solution; second, from (\ref{closed from a})-(\ref{closed from b}), we know that the homotopy algorithm needs to solve two linear systems of size $|J^{i}|$ at each step. Because the number of  support vectors is often small, that is, $\alpha^*$ is sparse, the homotopy algorithm solves smaller scale linear systems than IPM. This good property of $\alpha^*$ makes AL-Hom perform well in solving SVM optimizations.

\subsection{Contact problems of elasticity}
In this section, we solve the contact problems of elasticity  used as a benchmark in \cite{dostal2000solution, dostal2000duality,dostal2003augmented}
\begin{equation}\label{contact}
\begin{array}{l}
{\rm Minimize}~~~ q(u_1,u_2)=\sum_{i=1}^{2}\left(\int_{\Omega_i}|\nabla u_i|^2d\Omega -\int_{\Omega_i}Pu_i d\Omega \right),\\
{\rm subject~to}~~~ u_1(0,y)=0 ~{\rm and}~ u_1(1,y)\leq u_2(1,y) ~{\rm for } ~y\in[0,1],
\end{array}
\end{equation}
where $\Omega_1=(0,1)\times(0,1)$, $\Omega_2=(1,2)\times(0,1)$, $P(x,y)=-5$ for $(x,y)\in(0,1)\times[0.75,1)$, $P(x,y)=0$ for $(x,y)\in(0,1)\times(0,0.75)$,
$P(x,y)=-1$ for $(x,y)\in(1,2)\times(0,0.25)$ and $P(x,y)=0$ for $(x,y)\in(1,2)\times(0.25,1)$.

We followed Dost{\'a}l et al.  using finite difference to discretize (\ref{contact}) by  regular grids that are defined by the step size $h\in\{\frac{1}{32},\frac{1}{64},\frac{1}{128},\frac{1}{256},\frac{1}{512},\frac{1}{1024}\}$ in each direction in each subdomain $\Omega_i$.

The  discrete problem is a QP problem
\begin{equation}\label{disqp}
\begin{array}{l}
\min ~~ \frac{1}{2}x^{T}Qx+r^{T}x \\	
\rm{s.t.}~~~\emph{Ax}\leq 0,
\end{array}
\end{equation}
which is  transformed to form (\ref{scqp}) by introducing a slack variable $s$
\begin{equation}\label{disqp2}
\begin{array}{l}
\min ~~ \frac{1}{2}x^{T}Qx+r^{T}x \\	
\rm{s.t.}~~~\emph{Ax}+s=0\\
~~~~~~~s\geq 0.
\end{array}
\end{equation}
Thus, when $h=\frac{1}{1024}$, (\ref{disqp2}) has 2,101,250 variables and 1,025 equality constraints. We used PAL-Hom to solve (\ref{disqp2}) and compare it with the IPM solver in CPLEX. Moreover,  to show that the  strategy which uses the homotopy method to obtain the exact solutions of the augmented Lagrangian subproblems is valid, we also exactly asymptotically solved the augmented Lagrangian subproblems by the APG method. For convenience, we use PAL-APG to denote the augmented Lagrangian iterations with the subproblems exactly asymptotically solved by APG.
\vspace{-0.3cm}
\begin{table}[!htb]
\caption{\scriptsize{Experiments on the discrete contact problems of elasticity, ``Iter" denotes the number of the augmented Lagrangian subproblems solved. ``Total" denotes the computation time (seconds) of APG and the homotopy tracking steps together (for all augmented Lagrangian subproblems), ``Hom-tra." denotes the computation
time of the homotopy tracking steps. ``OM" denotes out of memory. The bolded computation times of PAL-Hom denotes that they are smaller than those of the other solvers.}}

\begin{center} \tiny
\begin{tabular}{llllllllllll}
\toprule
\multirow {2}{*}{$h$}&
\multirow {2}{*}{$m$}&
\multirow {2}{*}{$n$}&
\multicolumn{3}{c}{PAL-Hom} &
&
\multicolumn{2}{c}{PAL-APG}
&&
\multirow{2}{*}{IPM(cplex)}

  \\
\cline{4-6}
\cline{8-9}

   &&&Iter &Total&Hom-tra. & &Iter& Time&&  \\
\midrule
$1/32$ &33&  2,178&5&\bf{0.09}&0.02&&11&1.33&&0.32 \\
$1/64$ &65&  8,450&5&\bf{0.32}&0.06&&11&3.34&&1.34 \\
$1/128$ &129& 33,282&5&\bf{1.09}&0.18&&12&17.45&&13.33 \\
$1/256$  &257&   132,098&6&\bf{12.28}&0.85&&11&324.49&&127.42 \\
$1/512$&  513 &526,338&6&\bf{129.14}&10.09&&11&1651.90&&1408.52 \\
$1/1024$&  1,025 &2,101,250&6&\bf{873.17}&67.33&&12&12937.68&&OM \\
\bottomrule

\end{tabular}
\end{center}
\vspace{-0.6cm}
\label{table9}
\end{table}

From the results, we see that PAL-Hom requires fewer iterations than PAL-APG. Moreover the computation time of PAL-Hom is substantially smaller than that of PAL-APG. Furthermore, when the APG iteration obtains a low-precision solution, a good prediction of the optimal active set is obtained; therefore, the homotopy algorithm requires a small number of steps and little time to obtain the exact solution from the approximate solution. However, APG must continue iterating for the required precision, which requires much more time. The results demonstrate that   exactly solving the subproblems at the mid to end stage  by the homotopy algorithm is actually  valid and that PAL-Hom is substantially more efficient than IPM for solving this problem.

%%%%%%%%%%%%%%%%%%%%%%%%%%%%%%%%%%%%%%%%%%%%%%%%%%%%%%%%%%%%%%%%%%%%%%%%%%%%%%%%
\subsection{Randomly generated LPs and LPs from Netlib test set}
In this section, we solve LPs by PP-AL-Hom. We first randomly generated  LPs with MATLAB codes as follows.

~~~~~A=sprandn($m,n,d_A$); b=10$*$randn($m$,1); c=rand($n$,1).

%%%%%%%%%%%%%%%%%%%%%%%%%%%%%%%%%%%%%%%%%%%%%%%%%%%%%%%%%%%%%%%%%%%%%%%%%%%%%%%%%%%%%%
\begin{table}[!htb]

\caption{\scriptsize{Randomly generated LP}}
 \tiny
\begin{tabular}{p{2.6cm}p{2.5cm}p{2.5cm}p{2.5cm}p{2.5cm}p{2.5cm}p{2.5cm}}
\toprule
{Problem} &$m$& $n$&$d_A $\\
\hline
  LP-D1  &400&1000&1\\
  LP-D2  &800&2000&1\\
  LP-D3  &1000&5000&1\\
  LP-D4  &300&8000&1\\
  LP-D5  &4000&10000&1\\
  LP-S1  &100&2000&0.01\\
  LP-S2  &1000&5000&0.01\\
  LP-S3  &800&8000&0.01\\
  LP-S4  &4000&10000&0.01\\
  LP-S5  &800&15000&0.01\\
  LP-S6  &8000&20000&0.001\\
  LP-S7  &15000&32000&0.001\\
\toprule
\end{tabular}

\label{table9}
\end{table}
Additionally, we chose LPs from the Netlib test set. The chosen LPs have finite solutions and were up to a size $16558\times 49932$. For  randomly generated LPs,  PAL-Hom started from the original point, and for LPs from the Netlib test set, we used a projected Newton barrier method \cite{gill1986projected} to obtain an  approximate solution as an initial point, which would reduce the number of the  iterations (\ref{pg subproblem}).

\begin{table}[!htb]

\caption{\scriptsize{Experiments on randomly generated LP.  $f_{*}$, $f_{S}$ denote the optimal values obtained by PAL-Hom and the other corresponding solvers, ``OT'' denote the computation time (seconds) more than 25000s. The bolded computation times of PAL-Hom denotes  that they are smaller than those of the other solvers.}}
 \tiny

\begin{tabular}{p{0.7cm}p{1cm}p{1cm}p{1cm}p{1cm}p{1cm}p{1.6cm}}
\toprule
\multirow{1}{*}{Problem} &
Results&
\multicolumn{1}{l}{PP-AL-Hom} &
\multicolumn{1}{l}{IPM(cplex)}&
\multicolumn{1}{l}{Simplex(gurobi) } &
\multicolumn{1}{l}{IPM(matlab)}   &
\multicolumn{1}{l}{Simplex(matlab)} \\
   \midrule
   \multirow{3}{*}{LP-D1 } &Time &0.69 &0.80&0.58&3.33&17.90 \\
                           &$\| Ax-b\|$&7.8E-11 &4.0E-12& 4.8E-12& 1.3E-09&6.3E-10\\
                           &$f_{*}-f_{S}$ &-&-2.1E-11&-2.0E-11& -1.3E-10&-1.2E-11   \\
                           \midrule

   \multirow{3}{*}{LP-D2 } &Time &9.13 &8.47&5.59&30.25&232.99 \\
                           &$\| Ax-b\|$&2.4E-11&1.2E-11&1.3E-11&8.8E-13&3.8E-09 \\
                           &$f_{*}-f_{S}$ &-&-5.3E-12&-6.1E-12&-8.2E-12&1.0E-10  \\
                           \midrule

   \multirow{3}{*}{LP-D3 } &Time & 21.90&33.11&16.58&110.67&918.04 \\
                           &$\| Ax-b\|$&8.9E-10&1.2E-11&1.5E-11&9.9E-13&2.5E-09\\
                           &$f_{*}-f_{S}$ &-&1.0E-12&1.1E-12&9.9E-13&-1.0E-11  \\
                           \midrule

   \multirow{3}{*}{LP-D4 } &Time &11.36&3.82&2.27&19.47&55.78 \\
                           &$\| Ax-b\|$&5.7E-10&1.8E-12&1.6E-12&7.5E-11&3.8E-09\\
                           &$f_{*}-f_{S}$ &-&-6.2E-13&-6.2E-12&-5.0E-12&6.1E-12  \\
                           \midrule

   \multirow{3}{*}{LP-D5 } &Time &\bf{343.23}&1277.07&460.83&3472.62&OT \\
                           &$\| Ax-b\|$&2.2E-10&2.3E-10&2.6E-10&4.0E-12&-  \\
                           &$f_{*}-f_{S}$ &-&-5.0E-12&-1.2E-11&-7.7E-12&-  \\
                           \midrule

   \multirow{3}{*}{LP-S1 } &Time &0.09&0.00&0.00&0.02&0.07 \\
                           &$\| Ax-b\|$&2.1E-11&1.3E-13&9.7E-14&7.8E-13&4.9E-13\\
                           &$f_{*}-f_{S}$ &-&-2.1E-14&-2.1E-14&-2.1E-14&-2.3E-14  \\
                           \midrule

   \multirow{3}{*}{LP-S2 } &Time &2.43&1.62&0.93&6.40&350.84 \\
                           &$\| Ax-b\|$&1.0E-08&1.5E-10&4.5E-11&1.6E-08&2.7E-11\\
                           &$f_{*}-f_{S}$ &-&-8.1E-10&-8.1E-10&-1.4E-09&-5.3E-09 \\
                           \midrule

   \multirow{3}{*}{LP-S3 } &Time &5.16&1.09&0.39&3.37&171.34 \\
                           &$\| Ax-b\|$&8.0E-10&2.6E-10&6.8E-12&1.4E-10&2.7E-12\\
                           &$f_{*}-f_{S}$ &-&-5.3E-09&-5.3E09&-5.4E-10&-5.3E-10  \\
                           \midrule

   \multirow{3}{*}{LP-S4 } &Time & \bf{49.37}&96.59&59.78&458.33&OT \\
                           &$\| Ax-b\|$&4.2E-10&3.8E-09&2.9E-10&6.2E-09&-  \\
                           &$f_{*}-f_{S}$ &-&-5.2E-10&-4.9E-10&1.5E-10&-  \\
                           \midrule

   \multirow{3}{*}{LP-S5 } &Time &5.13&0.48&0.40&4.33&225.65 \\
                           &$\| Ax-b\|$&6.9E-11&6.2E-11&7.4E-12&4.2E-11&1.9E-10\\
                           &$f_{*}-f_{S}$ &-&4.2E-12&1.7E-13&3.7E-13&-4.3E-13  \\
                           \midrule

   \multirow{3}{*}{LP-S6 } &Time &311.45&148.93&97.58&3007.46&OT \\
                           &$\| Ax-b\|$&6.2E-09&2.7E-08&5.5E-09&1.4E-11&-  \\
                           &$f_{*}-f_{S}$ &-&-2.8E-08&-7.7E-09&-1.8E-09&-  \\
                           \midrule

   \multirow{3}{*}{LP-S7 } &Time &\bf{903.11}&2024.15&1289.04&24726.41&OT \\
                           &$\| Ax-b\|$&2.6E-09&2.1E-07&1.6E-08&1.6E-11&-  \\
                           &$f_{*}-f_{S}$ &-&-1.1E-06&-1.2E-06&-1.3E-06&- \\
\toprule
\end{tabular}

\label{table10}

\end{table}

%%%%%%%%%%%%%%%%%%%%%%%%%%%%%%%%%%%%%%%%%%%%%%%%%%%%%%%%%%%%%%%%%%%%%%%%%%%%%%%%%%%%%%%%%%%%

%%%%%%%%%%%%%%%%%%%%%%%%%%%%%%%%%%%%%%%%%%%%%%%%%%%%%%%%%%%%%%%%%%%%%%%%%%%%%%%%%%%%%%%%%%%%%%%%%
\begin{table}[!htb]
\caption{\scriptsize{Experiments on LPs from Netlib test set: part I (seconds). $f_{*}$, $f_{S}$ denote the optimal values obtained by PAL-Hom and the other corresponding solvers.}}
 \tiny
 \begin{tabular}{p{1.0cm}p{0.55cm}p{0.5cm}p{1cm}p{1cm}p{1cm}p{1cm}p{1cm}p{1.4cm}}
\toprule
\multirow{2}{*}{Problem} &
\multirow{2}{*}{m}&
\multirow{2}{*}{n}&
\multirow{2}{*}{Results}&
\multirow{2}{*}{PP-AL-Hom} &
\multicolumn{1}{l} {cplex}&
\multicolumn{1}{l}{gurobi} &
\multicolumn{2}{c}{matlab} \\
\cline{8-9}
&&&&&IPM&Simplex&IPM&Simplex\\

\midrule
   \multirow{3}{*}{adlittle}&
   \multirow{3}{*}{57 }&
   \multirow{3}{*}{138 }&Time &0.39&0.01&0.00&0.01&0.04 \\
                           &&&$\| Ax-b\|$&8.8E-11&2.4E-13&8.1E-14&3.1E-11&3.7E-13\\
                           &&&$f_{*}-f_{S}$&-&-1.4E-08&-1.2E-08&-1.2E-08&-1.2E-08\\
\midrule
   \multirow{3}{*}{afiro }&
   \multirow{3}{*}{27 }&
   \multirow{3}{*}{51 }&Time  & 0.04&0.00&0.00&0.01&0.01 \\
                           &&&$\| Ax-b\|$& 1.4E-13&1.4E-14& 1.4E-14&1.5E-12&1.1E-13\\
                           &&&$f_{*}-f_{S}$&-&1.1E-13& 01.1E-13&0.0E+00&1.1E-13\\

\midrule
   \multirow{3}{*}{agg2}&
   \multirow{3}{*}{516 }&
   \multirow{3}{*}{758 }&Time &0.51&0.01&0.01&0.08&1.27 \\
                           &&&$\| Ax-b\|$&1.8E-08&1.3E-13&3.4E-10&5.4E-10&1.3E-10\\
                           &&&$f_{*}-f_{S}$&-&-6.9E-04&-6.9E-04&-6.9E-04&-6.9E-04\\

\midrule
   \multirow{3}{*}{beaconfd }&
   \multirow{3}{*}{173 }&
   \multirow{3}{*}{295 }&Time  &0.27&0.00&0.00&0.02&0.02 \\
                           &&&$\| Ax-b\|$&4.3E-09&4.1E-11&1.2E-11&5.1E-11&3.2E-11 \\
                           &&&$f_{*}-f_{S}$&-&2.3E-04&2.3E-04&2.3E-04&2.3E-04 \\
\midrule
   \multirow{3}{*}{blend }&
   \multirow{3}{*}{74 }&
   \multirow{3}{*}{114 }&Time  &0.09&0.00&0.00&0.01&0.03 \\
                           &&&$\| Ax-b\|$&2.2E-09&63.9E-14&3.7E-13&6.0E-12&3.9E-14 \\
                           &&&$f_{*}-f_{S}$&-&-3.2E-07&-3.2E-07&-3.2E-07&-3.2E-07 \\

\midrule
   \multirow{3}{*}{d6cube }&
   \multirow{3}{*}{415 }&
   \multirow{3}{*}{6184 }&Time &49.13&0.08&0.07&0.40&16.93 \\
                           &&&$\| Ax-b\|$&8.6E-09&4.3E-11&7.8E-12&1.9E-09&5.1E-11 \\
                           &&&$f_{*}-f_{S}$&-&4.22E-07&4.22E-07&4.22E-07&4.22E-07 \\

\midrule
   \multirow{3}{*}{degen2 }&
   \multirow{3}{*}{444 }&
   \multirow{3}{*}{754 }&Time  &0.68&0.02&0.02&0.04& 2.04 \\
                           &&&$\| Ax-b\|$&1.6E-09&3.6E-15&3.6E-15&1.2E-12&2.6E-14 \\
                           &&&$f_{*}-f_{S}$&-3.7E-06&-3.7E-06&-3.7E-06&-3.7E-06&-3.7E-06 \\

\midrule
   \multirow{3}{*}{degen3 }&
   \multirow{3}{*}{1503 }&
   \multirow{3}{*}{2604 }&Time  & 17.19&0.31&0.10&0.79& 59.33 \\
                           &&&$\| Ax-b\|$&6.8E-09&1.1E-14&1.5E-14&7.2E-09&1.1E-13 \\
                           &&&$f_{*}-f_{S}$&-&-1.88E-05&-1.88E-05&-1.88E-05&-1.88E-05 \\

\midrule
   \multirow{3}{*}{maros-r7 }&
   \multirow{3}{*}{3136 }&
   \multirow{3}{*}{9408 }&Time &2.85&0.46&0.25& 3.99& 51.18 \\
                           &&&$\| Ax-b\|$&8.7E-09&4.6E-09&4.7E-09&1.8E-10&7.7E-08 \\
                           &&&$f_{*}-f_{S}$&-&-8.1E-08&-8.2E-08&-8.3E-08&-8.4E-08 \\
\midrule
   \multirow{3}{*}{psd$\_$02 }&
   \multirow{3}{*}{2953 }&
   \multirow{3}{*}{7716 }&Time  &5.91&0.03&0.02&0.18&2.15 \\
                           &&&$\| Ax-b\|$&7.4E-10&0.0E+00&0.0E+00&0.0E+00&0.0E+00 \\
                           &&&$f_{*}-f_{S}$&-&0.0E+00&0.0E+00&0.0E+00&0.0E+00 \\

\midrule
   \multirow{3}{*}{psd$\_$06 }&
   \multirow{3}{*}{9881 }&
   \multirow{3}{*}{29351 }&Time  & 26.44&0.16&0.12&6.41& 17.81 \\
                           &&&$\| Ax-b\|$&6.7E-09&0.0E+00&0.0E+00&0.0E+00&0.0E+00 \\
                           &&&$f_{*}-f_{S}$&-&-9.0E-04&-9.0E-04&-9.0E-04& -9.0E-04\\

\midrule
   \multirow{3}{*}{psd$\_$10 }&
   \multirow{3}{*}{16558}&
   \multirow{3}{*}{49932 }&Time  &213.32&0.39&0.21&31.66&50.83 \\
                           &&&$\| Ax-b\|$&2.0E-10&0.0E+00&0.0E+00&0.0E+00&0.0E+00 \\
                           &&&$f_{*}-f_{S}$&-&-2.8E-03&-2.8E-03&-2.8E-03&-2.8E-03 \\

\midrule
   \multirow{3}{*}{qap8 }&
   \multirow{3}{*}{912 }&
   \multirow{3}{*}{1632 }&Time &1.16&0.22&0.46&0.73&15.11 \\
                           &&&$\| Ax-b\|$&1.3E-09&5.0E-13&1.7E-12&1.1E-14&2.3E-14 \\
                           &&&$f_{*}-f_{S}$&-&3.3E-09&3.3E-09&3.7E-09&9.3E-09 \\

\midrule
   \multirow{3}{*}{qap12 }&
   \multirow{3}{*}{3192 }&
   \multirow{3}{*}{8856 }&Time  & 77.71&1.64&1.10&1506.63& 1342.79 \\
                           &&&$\| Ax-b\|$&8.6E-10&2.6E-12&1.9E-12&7.9E-09&1.9E-12 \\
                           &&&$f_{*}-f_{S}$&-&4.1E-06&4.1E-06&3.9E-06&4.1E-06 \\

\midrule
   \multirow{3}{*}{scorpion  }&
   \multirow{3}{*}{388 }&
   \multirow{3}{*}{466 }&Time  &0.54&0.01&0.01&0.02&0.22 \\
                           &&&$\| Ax-b\|$&3.1E-13&1.2E-15&8.1E-16&2.6E-15&1.2E-15 \\
                           &&&$f_{*}-f_{S}$&-&-5.8E-05&-5.8E-05&-5.8E-05&-5.8E-05 \\
\midrule
   \multirow{3}{*}{scsd1 }&
   \multirow{3}{*}{77 }&
   \multirow{3}{*}{760 }&Time  &0.32&0.01&0.01&0.01&0.09 \\
                           &&&$\| Ax-b\|$&9.1E-13&1.5E-16&1.1e-16&2.2e-13&2.1e-16 \\
                           &&&$f_{*}-f_{S}$&-&1.6E-11&1.6E-11&-2.5E-10&1.6E-11 \\
\midrule
   \multirow{3}{*}{scsd6 }&
   \multirow{3}{*}{147 }&
   \multirow{3}{*}{1350 }&Time  &0.29&0.01&0.02&0.02& 0.28 \\
                                &&&$\| Ax-b\|$&1.9E-12&5.9E-16&3.7E-16&1.9E-13&6.0E-16 \\
                                &&&$f_{*}-f_{S}$&-&-1.1E-09&-1.1E-09&-9.0E-09&2.4E-09 \\

\midrule
   \multirow{3}{*}{scsd8  }&
   \multirow{3}{*}{397 }&
   \multirow{3}{*}{2750 }&Time  & 0.19&0.02&0.04&0.02& 0.78 \\
                        &&&$\| Ax-b\|$&5.9E-11&3.2E-14&3.1E-13&3.1E-13&4.2E-14 \\
                           &&&$f_{*}-f_{S}$&-&-1.2E-08&1.5E-08&1.5E-08&1.5E-08 \\

\midrule
   \multirow{3}{*}{sctap1 }&
   \multirow{3}{*}{300 }&
   \multirow{3}{*}{660 }&Time &2.31&0.01&0.01&0.03&0.34 \\
                           &&&$\| Ax-b\|$&4.7E-10&5.5E-15&2.5E-15&6.0E-11&1.2E-11 \\
                           &&&$f_{*}-f_{S}$&-&1.9E-07&1.9E-07&1.9E-07&1.9E-07 \\

\midrule
   \multirow{3}{*}{sctap2 }&
   \multirow{3}{*}{1090 }&
   \multirow{3}{*}{2500 }&Time    &2.13&0.01&0.02&0.06& 5.10 \\
                           &&&$\| Ax-b\|$&7.3E-10&1.6E-14&8.9E-16&1.6E-12&3.3E-13 \\
                           &&&$f_{*}-f_{S}$&-&-2.5E-07&-2.5E-07&-2.5E-07&-2.5E-07 \\
\midrule
   \multirow{3}{*}{sctap3 }&
   \multirow{3}{*}{1480 }&
   \multirow{3}{*}{3340 }&Time  & 2.32&0.03&0.02&0.06& 6.86 \\
                           &&&$\| Ax-b\|$&1.2E-09&8.9E-15&6.2E-15&3.2E-12&3.2E-13 \\
                           &&&$f_{*}-f_{S}$&-&-6.5E-07&-6.5E-07&-6.5E-07& -6.5E-07\\

\toprule
\end{tabular}

\label{table11}
\end{table}
%%%%%%%%%%%%%%%%%%%%%%%%%%%%%%%%%%%%%%%%%%%%%%%%%%%%%%%%%%%%%%%%%%%%%%%
%%%%%%%%%%%%%%%%%%%%%%%%%%%%%%%%%%%%%%%%%%%%%%%%%%%%%%%%%%%%%%%%%%%%%%%%%%%%%%%%%%%%%%%%%%%%%%%%%
\begin{table}[!htb]
\caption{\scriptsize{Experiments on LPs from Netlib test set: part II (seconds). $f_{*}$, $f_{S}$ denote the optimal values obtained by PAL-Hom and the other corresponding solvers.}}

 \tiny
\begin{tabular}{p{1.0cm}p{0.55cm}p{0.5cm}p{1cm}p{1cm}p{1cm}p{1cm}p{1cm}p{1.0cm}}
\toprule
\multirow{2}{*}{Problem} &
\multirow{2}{*}{m}&
\multirow{2}{*}{n}&
\multirow{2}{*}{Results}&
\multirow{2}{*}{PP-AL-Hom} &
\multicolumn{1}{l} {cplex}&
\multicolumn{1}{l}{gurobi} &
\multicolumn{2}{c}{matlab} \\
\cline{8-9}
&&&&&{\hspace{0.1cm}IPM}&Simplex&IPM&Simplex\\

 \midrule
   \multirow{3}{*}{ship04l }&
   \multirow{3}{*}{402 }&
   \multirow{3}{*}{2166 }&Time  &0.41&0.01&0.01&0.02&0.25 \\
                           &&&$\| Ax-b\|$&3.5E-10&4.4E-14&4.9E-13&4.4E-11&2.3E-14 \\
                           &&&$f_{*}-f_{S}$&-&-7.5E-05&-3.1E-05&-3.1E-05&-3.1E-05 \\

\midrule
   \multirow{3}{*}{ship04s }&
   \multirow{3}{*}{402 }&
   \multirow{3}{*}{1506 }&Time &0.30&0.01&0.01&0.02&0.09 \\
                           &&&$\| Ax-b\|$&1.7E-09&7.7E-14&2.9E-14&6.8E-09&6.6E-14 \\
                           &&&$f_{*}-f_{S}$&-&-4.2E-04&-4.2E-04&-4.2E-04&-4.2E-04 \\
\midrule
   \multirow{3}{*}{ship08l }&
   \multirow{3}{*}{778 }&
   \multirow{3}{*}{4363 }&Time &3.38&0.01&0.01&0.05&0.45 \\
                           &&&$\| Ax-b\|$&2.3E-12&4.7E-14&3.2E-14&2.2E-10&1.7E-13 \\
                           &&&$f_{*}-f_{S}$&-&-6.5E-07&-1.4E-07&-1.4E-07&-1.4E-07 \\
\midrule
   \multirow{3}{*}{ship08s }&
   \multirow{3}{*}{778 }&
   \multirow{3}{*}{2476 }&Time &2.13&0.02&0.01&0.03&0.18 \\
                           &&&$\| Ax-b\|$&1.9E-12&2.8E-14&1.8E-14&2.8E-11&1.0E-11 \\
                           &&&$f_{*}-f_{S}$&-3.6E-08&1.1E-07&1.1E-07&1.1E-07&1.1E-07 \\

\midrule
   \multirow{3}{*}{ship12l}&
   \multirow{3}{*}{1151 }&
   \multirow{3}{*}{5533 }&Time &6.33&0.02&0.02&0.06&0.67 \\
                           &&&$\| Ax-b\|$&4.4E-12&3.6E-14&3.8E-14&3.3E-11&3.6E-13 \\
                           &&&$f_{*}-f_{S}$&-&-2.2E-07&-1.8E-07&-1.8E-07&-1.8E-07 \\
\midrule
   \multirow{3}{*}{mship12s }&
   \multirow{3}{*}{1151 }&
   \multirow{3}{*}{2869 }&Time &2.64&0.01&0.02&0.02&0.28 \\
                           &&&$\| Ax-b\|$&1.8E-11&4.9E-13&6.3E-14&3.2E-11&1.1E-13 \\
                           &&&$f_{*}-f_{S}$&-&-1.3E-05&-9.5E-07&-9.5E-07&-9.5E-07 \\
\midrule
   \multirow{3}{*}{truss}&
   \multirow{3}{*}{1000 }&
   \multirow{3}{*}{8806 }&Time& 2.67&0.07&1.91&0.18& 21.66 \\
                           &&&$\| Ax-b\|$&1.7E-09&2.1E-13&1.9E-13&1.8E-11&1.0E-11 \\
                           &&&$f_{*}-f_{S}$&-&7.2E-06&7.2E-06&7.2E-06&7.2E-06 \\

\toprule
\end{tabular}

\label{table12}
\end{table}
%%%%%%%%%%%%%%%%%%%%%%%%%%%%%%%%%%%%%%%%%%%%%%%%%%%%%%%%%%%%%%%%%%%%%%%
We report the results in Tables \ref{table10}-\ref{table12} and the time of PAL-Hom in Tables \ref{table11}-\ref{table12} has included the computation time of the projected Newton barrier method. The results show that PAL-Hom is able to solve the randomly generated LPs and LPs from the Netlib test set. For randomly generated LPs, PAL-Hom is competitive with the other solvers. For LPs from the Netlib test set, PAL-Hom is not as good as the IPM solvers in CPLEX and MATLAB, and the simplex solver in Gurobi, but for some problems, PAL-Hom is more effective than the simplex solver in MATLAB.

%%%%%%%%%%%%%%%%%%%%%%%%%%%%%%%%%%%%%%%%%%%%%%%%%%%%%%%%%%%%%%%%%%%%%%%%%%%%%%%%%%
\section{Conclusion}
In this paper, we present a PAL-Hom (AL-Hom) algorithm for convex QP problems, which takes the proximal ALM as the outer iteration and the homotopy algorithm as the inner iteration. Compared with IPM, AS and PAS, the size of the KKT systems  solved in PAL-Hom is much smaller, especially when the solution is sparse such as in the problems from SVM. Moreover, compared with PAS, the KKT systems in the tracking steps of PAL-Hom would always be invertible so that we do not need to exchange indices to keep the invertibility as in qpOASES. Furthermore, it is substantially easier to design an efficient warm start for PAL-Hom than  for the QP problem (\ref{scqp}) (PAS). Although  we do not pay significant attention to optimizing the codes, PAL-Hom is shown to be faster than the IPM solver in CPLEX for certain problems, such as randomly generated QPs, LPs and some QPs in the CUTEr test. In particular, SVM  problems and the discrete contact problems of elasticity, PAL-Hom is more than 10 times faster than IPM. Given this   practical performance, we believe that our algorithms are promising.

The presented homotopy algorithm is shown to be efficient for nonnegative QP problems (augmented Lagrangian subproblems) for the following reasons. First, APG is effective at predicting the optimal active set, which provides a good warm start for the homotopy algorithm. With the warm start, the homotopy algorithm often needs fewer iterations to obtain an exact solution.  The Cholesky factor update technique improves the performance of the homotopy algorithm by reducing the computation of solving the KKT systems. Moreover, benefiting from  the $\varepsilon$-precision verification and correction technique that address the incorrect update of the active set caused by   large condition numbers and a lack of strict complementarity, the homotopy algorithm is shown to be robust for the augmented Lagrangian problems with large condition numbers. The numerical results demonstrate that the homotopy algorithm is substantially more efficient than PAS, ASA and IPM in solving the augmented Lagrangian subproblems.

Simultaneously, based on the AL-Hom method, we use PP-AL-Hom to solve the LP which is proved to converge in a finite number of steps. Moreover, the estimate of the number of  maximum iterations and the descent of the objective are presented. The numerical results show that PP-AL-Hom is competitive to IPM in solving randomly generated problems.

\section*{Acknowledgments}

The authors would like to thank Xiaoliang Song (School of Mathematical Sciences, Dalian University of Technology) for his valuable suggestions, which led to improvement in this paper. This research was supported by the National Natural Science Foundation of China (11571061, 11401075 and 11701065) and the Fundamental Research Funds for the Central Universities (DUT16LK05 and DUT17LK14)

%\begin{acknowledgements}
%\end{acknowledgements}

\bibliographystyle{spmpsci}      % mathematics and physical sciences
%\bibliography{lpreferences}
% etc
%\end{thebibliography}

\end{document}